\documentclass[preprint,authoryear,12pt]{elsarticle}

\usepackage{amssymb}

\usepackage{amsthm}
\usepackage[sumlimits]{amsmath}
\usepackage{amsfonts}
\usepackage{color,rotating}
\usepackage[margin=2cm]{geometry}
\usepackage{graphicx}
\usepackage{psfrag}
\psfrag{test}{$P1^3$}

\DeclareMathOperator{\diff}{d}

\newcommand{\iso}{{\mathrm{iso}}}
\newcommand{\PP}{{\mathrm{P}}}
\newtheorem{theorem}{Theorem}

\newtheorem{definition}[theorem]{Definition}
\newtheorem{example}[theorem]{Example}

\def\MM#1{\boldsymbol{#1}}
\newcommand{\pp}[2]{\frac{\partial #1}{\partial #2}} 

\newcommand{\dd}[2]{\frac{\diff#1}{\diff #2}}

\newcommand{\U}{\MM{u}}
\newcommand{\W}{\MM{w}}
\newcommand{\LL}{\left\langle}
\newcommand{\RR}{\right\rangle}

\usepackage{xspace}
\newcommand{\pdgp}{${P}1_{{DG}}$-${P}2$\xspace}
\newcommand{\pdg}{${P}1_{{DG}}$\xspace}
\newcommand{\ptwo}{${P}2$\xspace}
\newcommand{\hptwo}{$H({P}2)$\xspace}
\newcommand{\ptwobar}{$\overline{P2}$\xspace}
\newcommand{\pndgpn}{${P}n_{{DG}}$-${P}(n+1)$\xspace}

\journal{Journal of Computational Physics}

\begin{document}

\begin{frontmatter}



\title{Numerical wave propagation for the triangular \pdgp
 finite element pair} 


\author[aero]{C.~J. Cotter}
 \author[ese,grantham]{D.~A. Ham} 

\address[aero]{Department of Aeronautics, Imperial College London,
 South Kensington Campus, London SW7 2AZ}
\address[ese]{Department of Earth Science and Engineering, Imperial College London,
South Kensington Campus, London SW7 2AZ}
\address[grantham]{Grantham Insitute for Climate Change, Imperial College London,
South Kensington Campus, London SW7 2AZ}

\begin{abstract}
  The $f$-plane and $\beta$-plane wave propagation properties are
  examined for discretisations of the linearized rotating
  shallow-water equations using the \pdgp finite element pair on
  arbitrary triangulations in planar geometry. A discrete Helmholtz
  decomposition of the functions in the velocity space based on
  potentials taken from the pressure space is used to provide a
  complete description of the numerical wave propagation for the
  discretised equations. In the $f$-plane (planar geometry, Coriolis
  force independent of space) case, this decomposition is used to
  obtain decoupled equations for the geostrophic modes, the
  inertia-gravity modes, and the inertial oscillations. As has been
  noticed previously, the geostrophic modes are steady. The Helmholtz
  decomposition is used to show that the resulting inertia-gravity
  wave equation is third-order accurate in space. In general the \pdgp
  finite element pair is second-order accurate, so this leads to very
  accurate wave propagation. It is further shown that the only
  spurious modes supported by this discretisation are spurious
  inertial oscillations which have frequency $f$, and which do not
  propagate. A restriction of the $P1_{DG}$ velocity space is proposed
  in which these modes are not present, leading to a finite element
  discretisation which is completely free of spurious modes.  The
  Helmholtz decomposition also allows a simple derivation of the
  quasi-geostrophic limit of the discretised \pdgp equations in the
  $\beta$-plane (planar geometry, Coriolis force linear in space) case
  resulting in a Rossby wave equation which is also third-order
  accurate. This means that the dispersion relation for the wave
  propagation is very accurate; an illustration of this is provided by
  a numerical dispersion analysis in the case of a triangulation
  consisting of equilateral triangles.
\end{abstract}


\end{frontmatter}

\section{Introduction}

Recently there has been growing interest in developing more general
horizontal discretisation schemes for numerical weather prediction
(NWP) models with computational meshes constructed from triangles or
hexagons. There are two principal motivations for this.  Firstly,
geodesic grids (which are obtained by iterative refinement of an
icosahedron using triangles, sometimes transforming to the dual grid
which is a mesh of hexagons with exactly 12 pentagons located at the
vertices of the original icosahedron) provide similar grid cell areas
over the entire sphere, which has possible advantages for accurate
representation of wave propagation. Furthermore, geodesic grids also
avoid the very fine grid cells obtained near the North and South poles
on latitude-longitude grids, which lead to large Courant numbers, and
cause bottlenecks in communication between processors on parallel
systems.  This has led to a number of groups developing weather and
climate models which use geodesic grids
\citep{Ri+2000,Ma+2002,Sa+2008}. Secondly, triangles facilitate the
implementation of adaptive mesh refinement.  This allows nested
regional models within a global model, and further allows dynamic mesh
refinement in which the mesh resolution is locally modified in
response to the dynamics in the course of a forecast. The development
of new numerical schemes that correctly represent the qualitative
properties of wave propagation on these grids, and under adaptive-mesh
refinement, is crucial.

Possible discretisations on triangular or hexagonal meshes are
obtained using three different approaches: finite difference methods,
finite volume methods and finite element methods. To eliminate
spurious pressure modes, finite difference methods use a C-grid in
which the edge-normal velocity is stored at the edge-centres, and the
pressure is stored at the cell-centres. On quadrilateral grids, the
wave propagation is observed to be well represented provided that the
Rossby radius is well-resolved \citep{ArLa77,Fo-Ra96,Ra94}. On
triangular and hexagonal grids the problem lies in finding a scheme
for reconstructing the Coriolis force (which requires the tangential
velocity) from the normal velocity. Recently, a reconstruction scheme
was found which results in steady geostrophic modes for C-grid
discretisations on the regular hexagonal grid in the plane
\citep{Th08}. In the same paper it was shown that the resulting
discrete system on the $\beta$-plane has a spurious extra Rossby wave
branch, with very slow Eastward phase velocities. This reconstruction
was extended to arbitrarily structured C-grids in \citet{ThRiSkLe2009}.
The finite element method provides a great degree of flexibility in
the choices of discretisation for velocity and pressure. Amongst the
many finite element pairs that have been proposed for the rotating
shallow-water equations are the $\PP1_{NC}$-$\PP1$ and $\PP1$-$\iso$ $\PP2
$-$P1$ elements (investigated and compared to several other element
pairs in \citet{Ro_etal1998}), the RT0 elements (introduced in
\citet{RaTh1977} and proposed for the shallow-water equations in
\citet{WaCa1998}) and equal-order elements with stabilisation (also
proposed in \citet{WaCa1998}). 

In this paper we study the wave propagation properties of the recently
proposed \pdgp finite element discretisation. This discretisation uses
a mixed finite element pair with The \pdgp finite element
discretisation was introduced in \citet{CoHaPaRe2009}, and was designed
to accomodate the geostrophic balance relation between pressure and
velocity without introducing spurious pressure modes. This is achieved
by using a quadratic (\ptwo) continuous finite element basis for
pressure, and a linear discontinuous (\pdg) finite element basis for
velocity (hence the name \pdgp). The pressure polynomials are one
order higher than the velocity polynomials, which accomodates the
geostrophic balance relation since the pressure gradient and the
velocity are both linear within each element. Making the velocity
basis discontinuous increases the number of velocity degrees of
freedom so that there are no spurious pressure modes. The lack of
pressure modes was investigated numerically in \citet{CoHaPaRe2009} and
subsequently proved in \citet{CoHaPa2009}, where it was also shown that
this combination of spaces means that geostrophically balanced states
are exact steady states of the linear equations on arbitrary
unstructured meshes (this property can also be obtained for C-grid
finite difference methods as described in \citet{ThRiSkLe2009}, with
the added restriction that the meshes satisfy an orthogonality
property). In this paper we go further and produce a complete
description of the numerical wave propagation properties of \pdgp,
which is facilitated by the construction of a discrete Helmholtz
decomposition of the \pdg space. 

The rest of this paper is organised as follows. In Section
\ref{helmholtz}, we show that \pdgp has a discrete Helmholtz
decomposition. In Section \ref{f plane} we use this decomposition to
analyse the wave propagation on the $f$-plane. We show that there are
three types of modes: steady geostrophic modes, inertia-gravity modes,
and inertial oscillations (of which only one is a physical mode). We
show that the inertial oscillations do not propagate and can be
filtered out by solving two discretised elliptic equations. We also
show that the velocity may be eliminated to obtain a third-order
accurate inertia-gravity wave equation, and hence claim that the wave
propagation is very accurate on arbitrary unstructured meshes. In
Section \ref{beta plane}, we use the Helmholtz decomposition to 
analyse the Rossby wave propagation on the $\beta$-plane in the
quasi-geostrophic limit (following the approach of \citet{Th08}). We
obtain a third-order accurate Rossby wave equation, and hence claim
that the Rossby wave equation is also very accurate. Finally, in
section \ref{summary} we give a summary and outlook.

\section{Discrete Helmholtz decomposition for \pdgp}
\label{helmholtz}
In this section we show that the \pdgp finite element discretisation
has a discrete Helmholtz decomposition for \pdgp. We shall adopt the
notation that the $\delta$ superscript indicates a numerical
approximation in a finite element space; functions without subscripts
indicate continuous fields.  We start by stating two properties of
\pdgp which we shall use throughout.
\begin{definition}[Embedding conditions]
\label{embedding}
Let $V$ be the chosen vector space of finite element velocity fields (in
the case of \pdgp, $V$ is the space \pdg of velocity fields $\MM{u}^\delta$ that 
are linear in each triangular element, with no continuity constraints across
element boundaries), and let $H$ be the chosen vector space of finite
element pressure fields (in the case of \pdgp, $H$ is the space \ptwo
of pressure fields $h^\delta$ that are quadratic in each triangular
element and are constrained to be continuous across element boundaries).
\begin{enumerate}
\item The operator $\nabla$ defined by the pointwise gradient
\[
{q}^\delta(\MM{x}) = \nabla h^\delta(\MM{x})
\]
maps from $H$ into $V$.
\item The skew operator $\perp$ defined by the pointwise
formula
\[
\MM{q}^\delta(\MM{x}) = (\MM{u}^\delta(\MM{x}))^\perp = (-u^\delta_2,u^\delta_1)
\]
maps from $V$ into itself.
\end{enumerate}
\end{definition}
These are the only conditions that we use in the paper and hence any
properties extend to any other finite element pair that satisfies
these conditions (P0-P1 or \pndgpn with any $n>1$, for example).

These conditions are most definitely not satisfied by all possible
pairs $(V,H)$, as illustrated by the following examples.
\begin{example}[P1-P1]
  The finite element pair known as P1-P1 (which may be used for the
  shallow-water equations but requires stabilisation as described in
  \citet{WaCa1998}) is defined as follows:
\begin{itemize}
\item The mesh $\mathcal{M}$ is composed of triangular elements.
\item $H$ is the space of elementwise-linear functions $h^\delta$
  which are continuous across element boundaries.
\item $V$ is the space of vector fields $\MM{u}^\delta$ with both of
  the Cartesian components $(u^\delta,v^\delta)$ in $H$.
\end{itemize} 
Condition 1 of Definition \ref{embedding} is not satisfied by the
P1-P1 pair since gradients of functions in $H$ are discontinuous
across element boundaries. Condition 2 is satisfied since the same
continuity conditions are required for normal and tangential
components.
\end{example}
\begin{example}[RT0]
  The lowest order Raviart-Thomas \citep{RaTh1977} velocity space
  (known as RT0) is constructed on a mesh $\mathcal{M}$ composed of
  triangular elements. It consists of elementwise constant vector
  fields which are constrained to have continuous normal components
  across element boundaries. RT0 does not satisfy condition 2 of
  Definition \ref{embedding} since the $\perp$ operator transforms
  vector fields with discontinuities in the tangential component
  (which are permitted in RT0) into vector fields with discontinuities
  in the normal component (which are not).
\end{example}

We now describe some examples of finite element pairs which \emph{do}
satisfy the conditions in Definition \ref{embedding}.
\begin{example}[P0-P1]
  The finite element pair known as P0-P1 (applied to ocean modelling
  in \citet{Um+2004}, and analysed in \citet{LeRoPo2007,LePo2008}) is defined as follows:
\begin{itemize}
\item The mesh $\mathcal{M}$ is composed of triangular elements.
\item $H$ is the space of elementwise-linear functions $h^\delta$
  which are continuous across element boundaries.
\item $V$ is the space of elementwise-constant vectors with discontinuities
  across element boundaries permitted. 
\end{itemize}
\end{example}
\begin{example}[\pdgp]
  The finite element pair known as \pdgp \citep{CoHaPaRe2009} is
  defined as follows:
\begin{itemize}
\item The mesh $\mathcal{M}$ is composed of triangular elements.
\item $H$ is the space of elementwise-quadratic functions $h^\delta$
  which are continuous across element boundaries.
\item $V$ is the space of elementwise-linear vectors with discontinuities
  across element boundaries permitted. 
\end{itemize}
\end{example}
Each of these examples satisfy both conditions in Definition
\ref{embedding}: condition 1 holds because taking the gradient of a
elementwise polynomial $n-1$ which is continuous across element
boundaries results in a vector field which is discontinuous across
element boundaries and is composed of elementwise polynomials of one
degree $n$, and condition 2 holds since the velocity space uses the
same continuity constraints for normal and tangential components
\emph{e.g.}  both components are allowed to be discontinuous. This
defines a whole sequence of high-order \pndgpn element pairs. Similar
elements can be constructed on quadrilateral elements.  Since we only
require these two conditions to prove our optimal balance property
which holds on arbitrary meshes, we can also construct finite element
spaces on mixed meshes composed of quadrilaterals and triangles, for
example. It is also possible to use $p$-adaptivity in which different
orders of polynomials are used in different elements, as long as the
conditions are satisfied. To make the rest of the paper less abstract,
we shall only discuss \pdgp, but all of the results are easily
extended (with th appropriate orders of accuracy) to any element pair
satisfying Definition \ref{embedding}.

Next we note that the gradient and skew-gradient any two pressure
fields $\phi^\delta$, $\psi^\delta$ in the pressure space \ptwo are
orthogonal in the $L_2$ inner product,
\[
\langle \nabla\psi^\delta,\nabla^\perp \phi^\delta\rangle
=\int_{\Omega} \nabla\psi^\delta \cdot \nabla^\perp \phi^\delta
\diff{V} =0,
\]
where $\Omega$ is the solution domain which is either the sphere, or
periodic boundary conditions. This was proved by direct computation in
\citet{CoHaPa2009}.  Hence, any velocity field $\MM{u}^\delta$ in \pdg
can be written uniquely in an orthogonal decomposition
\begin{equation}
\label{discrete helmholtz decomposition}
\MM{u}^\delta = \bar{\MM{u}}^\delta + \nabla\phi^\delta +
\nabla^\perp\psi^\delta + \hat{\MM{u}}^\delta,
\end{equation}
where $\bar{\MM{u}}^\delta$ is independent of space,
where $\phi^\delta$ and $\psi^\delta$ are both in the space \ptwobar,
which consists of \ptwo functions with mean zero, \emph{i.e.}
\[
\LL \phi^\delta, 1 \RR = 
\int_{\Omega}\phi^\delta\diff{V} = 0,
\quad \LL \psi^\delta , 1 \RR = 
\int_{\Omega}\psi^\delta\diff{V}=0,
\]
and where $\hat{\MM{u}}^\delta$ is orthogonal to the gradient or
skew-gradient of any \ptwobar function $\alpha^\delta$, \emph{i.e.}
\[
\LL\hat{\MM{u}}^\delta,\nabla\alpha^\delta\RR =
\LL\hat{\MM{u}}^\delta,\nabla^\perp\alpha^\delta\RR =0.
\]
Furthermore, if any such $\hat{\MM{u}}^\delta$ satisfies
\[
\LL\hat{\MM{u}}^\delta,\hat{\MM{u}}^\delta\RR=0,
\]
then $\hat{\MM{u}}^\delta=\MM{0}$, since $\hat{\MM{u}}^\delta$ is
obtained from orthogonal completion. In general the dimension of the
orthogonal subspace containing the vector fields of the form
$\hat{\MM{u}}^\delta$ is non-zero, since there are more than twice as
many degrees of freedom in the velocity space $V$ as the pressure
space $H$. The dimension of $V$ is $6n_f$ (where $n_f$ is the number
of elements), and the dimension of $F$ is $n_v+n_e$ (where $n_v$ is
the number of vertices and $n_e$ is the number of edges). For doubly
periodic boundary conditions, Euler's polyhedral formula on the torus
then gives $\dim(H)=n_v+n_e=2n_e-n_f$. For a triangulation, $2n_e=3n_f$
since each triangle has three edges which are each shared between two
faces, so $\dim(H)=2n_f<3n_f=\dim(V)/2$. Since $2\dim(H)<\dim(V)$ it
is not possible to span $V$ entirely with functions of the form
$\nabla^\perp\psi^\delta+\nabla\phi^\delta$,
$\psi^\delta,\phi^\delta\in H$, and so components of the form
$\hat{\U}^\delta$ will always be present.

Equation \eqref{discrete helmholtz decomposition} is identical to the
Helmholtz decomposition for arbitrary continuous velocity fields in
which any continuous velocity field $\MM{u}$ can be written as a
constant plus a gradient of a potential plus the skew gradient of a
streamfunction; the only difference in the discrete \pdgp case is the
extra component $\hat{\MM{u}}^\delta$. This extra component gives rise to
spurious inertial oscillations in the \pdgp finite element
discretisation applied to the rotating shallow-water equations. It is
possible to describe a reduced velocity space, which we call \hptwo,
consisting of velocity fields which can be written as
\[
\MM{v}^\delta = \bar{\MM{v}}^\delta + \nabla\phi^\delta +
\nabla^\perp\psi^\delta,
\]
where $\bar{\MM{v}}^\delta$ is independent of space,
where $\phi^\delta$ and $\psi^\delta$ are both in the space \ptwobar,
\emph{i.e.} we have removed the spurious velocity component. It is
possible to project a \pdg velocity field $\MM{u}^\delta$ into \hptwo,
by first computing the mean component,
\[
\bar{\MM{u}}^\delta = \frac{\int_{\Omega}\MM{u}^\delta\diff{V}}
{\int_{\Omega}\diff{V}},
\]
and then extracting the velocity potential and streamfunction by
solving
\[
\LL\nabla\alpha^\delta,\nabla\phi^\delta\RR
= \LL\nabla\alpha^\delta,\MM{u}^\delta\RR,
\]
and
\[
\LL\nabla\alpha^\delta,\nabla\psi^\delta\RR
= \LL\nabla^\perp\alpha^\delta,\MM{u}^\delta\RR,
\]
for all \ptwobar functions $\alpha^\delta$. This amounts to solving 
elliptic problems for $\phi^\delta$ and $\psi^\delta$. Then, the
projection of $\MM{u}^\delta$ into \hptwo is given by 
\[
\bar{\MM{u}}^\delta + \nabla\phi^\delta + \nabla^\perp\psi^\delta.
\]

\section{Discrete wave propagation on the $f$-plane}
\label{f plane}
In this section we describe all of the numerical solutions obtained
from \pdgp applied to the $f$-plane.

\subsection{Discrete wave equation on the $f$-plane}

The \pdgp spatial discretisation of the rotating shallow-water
equations (see \citet{CoHaPa2009} for a derivation) is
\begin{eqnarray}
\label{u eqn}
\dd{}{t}\LL \MM{w}^\delta,\MM{u}^\delta\RR 
+ \LL f\MM{w}^\delta,(\MM{u}^\delta)^\perp\RR
& = &
-c^2\LL\MM{w}^\delta,\nabla \eta^\delta\RR, \\
\label{h eqn}
\dd{}{t}
\LL\phi^\delta,\eta^\delta \RR & = & 
\LL\nabla\phi^\delta,\MM{u}^\delta\RR, 
\end{eqnarray}
where the velocity $\MM{u}^\delta$ is in \pdg, the layer depth 
$\eta^\delta = H(1+\eta^\delta)$ is in \ptwo, for
all test functions $\MM{w}^\delta$ in \pdg and $\phi^\delta$ in \ptwo,
and where $c^2=gH$ is the non-rotating wave propagation speed, $g$ is
the acceleration due to gravity, $H$ is the mean layer depth and $f$
is the Coriolis parameter. 

On the $f$-plane, $f$ is a constant, and so we may take it outside the
Coriolis integral. Applying the discrete Helmholtz decomposition to
the velocity $\MM{u}^\delta$ and the velocity test functions
$\MM{w}^\delta$, \emph{i.e.},
\[
\U^{\delta} = \bar{\U}^\delta+ \nabla\phi^\delta + \nabla^\perp\psi^\delta + 
\hat{\U}^\delta, \quad
\W^{\delta} = \bar{\W}^\delta+ \nabla\alpha^\delta + \nabla^\perp\beta^\delta + 
\hat{\W}^\delta,
\]
equations (\ref{u eqn}-\ref{h eqn}) become (after removing products of
orthogonal quantities)
\begin{eqnarray}
\label{f plane phi}
\dd{}{t}\LL \nabla\alpha^\delta, \nabla \phi^\delta \RR
- f\LL \nabla\alpha^\delta, \nabla\psi^\delta \RR + c^2\LL \nabla\alpha^\delta, \nabla\eta^\delta \RR
& = & 0,  \\
\label{f plane psi}
\dd{}{t}\LL\nabla\alpha^\delta, \nabla\psi^\delta\RR + f\LL \nabla\alpha^\delta,
\nabla\phi^\delta \RR & = & 0, \\
\label{f plane eta}
\dd{}{t}\LL \alpha^\delta, \eta^\delta \RR
 - \LL\nabla\alpha^\delta,\nabla\phi^\delta\RR & = & 0, \\
\label{f plane inertial}
\dd{}{t}\LL \bar{\W}^\delta,\bar{\U}^\delta \RR + f\LL\bar{\W}^\delta,\bar{\U}^\delta \RR & = & 0,
\\
\label{f plane spurious}
\dd{}{t}\LL \hat{\W}^\delta,\hat{\U}^\delta \RR + f\LL\hat{\W}^\delta,(\hat{\U}^\delta)^\perp \RR & = & 0,
\end{eqnarray}
These solutions exhibit four types of orthogonal modes: geostrophic balance,
inertia gravity waves, the physical inertial oscillation, and spurious
inertial oscillations due to the presence of $\hat{\MM{u}}$. We shall now
describe these modes one by one.

\subsection{Geostrophic balance}
\label{geostrophic balance}
For the continuous equations before discretisation, geostrophically
balanced modes are obtained from non-zero steady solutions of the
equations. As shown in \citet{CoHaPa2009}, in the \pdgp discretisation
solutions which satisfy the geostrophic balance relation are also
exactly steady. To see this within the framework of this paper, assume
a steady state, then equations (\ref{f plane phi}-\ref{f plane
  spurious}) become
\begin{eqnarray}
\label{f plane phi steady}
- f\LL \nabla\alpha^\delta, \nabla\psi^\delta \RR + c^2\LL \nabla\alpha^\delta, \nabla\eta^\delta \RR
& = & 0,  \\
\label{f plane psi steady}
f\LL \nabla\alpha^\delta,
\nabla\phi^\delta \RR & = & 0, \\
\label{f plane eta steady}
- \LL\nabla\alpha^\delta,\nabla\phi^\delta\RR & = & 0, \\
\label{f plane inertial steady}
f\LL\bar{\W}^\delta,\bar{\U}^\delta \RR & = & 0,
\\
\label{f plane spurious steady}
f\LL\hat{\W}^\delta,(\hat{\U}^\delta)^\perp \RR & = & 0.
\end{eqnarray}
Equations \eqref{f plane psi steady} and \eqref{f plane eta steady}
both imply that $\phi^\delta=0$ since they are the usual continuous
finite element discretisations of the Laplace equation which has no
non-zero solutions because $\phi^\delta$ and $\alpha^\delta$ are both
restricted to \ptwobar. Similarly equations \eqref{f plane inertial
  steady} and \eqref{f plane spurious steady} imply that
$\bar{\U}^\delta=\hat{\U}^\delta=\MM{0}$.  Equation \eqref{f plane phi steady} is
the discrete geostrophic balance relation between $\psi^\delta$ and
$\eta^\delta$, and the Laplace operator can be inverted (since the
finite element discretisation of the Poisson equation has a unique
solution for solutions in \ptwobar) to obtain the \emph{pointwise}
geostrophic balance relation
\[
f\psi^\delta = c^2\eta^\delta,
\]
as noted in \citet{CoHaPa2009}. This means that \pdgp has an excellent
representation of geostrophic balance.

\subsection{Inertia gravity waves}
\label{inertia gravity waves}
The physical wave variables $\phi^\delta$,
$\psi^\delta$ and $\eta^\delta$ are uncoupled to the mean velocity
component $\bar{\U}^\delta$ and the spurious velocity component $\hat{\U}^\delta$.
To obtain the discrete inertia gravity wave equation, the time
derivative applied to equation \eqref{f plane phi} gives
\[
\dd{^2}{t^2}\LL \nabla\alpha^\delta, \nabla \phi^\delta \RR
- f
\dd{}{t}\LL \nabla\alpha^\delta, \nabla\psi^\delta \RR +
\dd{}{t} c^2\LL \nabla\alpha^\delta, \nabla\eta^\delta \RR
= 0 .
\]
Substitution of equations \eqref{f plane psi} and \eqref{f plane eta}
then give
\begin{equation}
\left(\dd{^2}{t^2}+f^2\right)\dd{}{t}\LL \alpha^\delta, \eta^\delta \RR
+\dd{}{t} c^2\LL \nabla\alpha^\delta, \nabla\eta^\delta \RR
= 0 .
\label{ig equation}
\end{equation}
This is the usual continuous finite element discretisation of the
inertia-gravity wave equation
\begin{equation}
\left(\pp{^2}{t^2}+f^2\right)\pp{}{t}\eta - c^2\nabla^2\pp{\eta}{t} =0.
\end{equation}
Since only \ptwo functions are present, the solution $\eta^\delta$ is
third-order accurate, as opposed to the second-order accuracy expected
with a first-order velocity discretisation.  This higher-than-expected
accuracy means that \pdgp has a very accurate representation of
inertia-gravity wave propagation. In particular, it should be expected
that the phase velocity is more independent of mesh orientation than
other second-order methods.  The equivalent property for $P0$-$P1$ was
noted in \citet{LeRoPo2007}, namely that the inertia-gravity dispersion
relation was one order more accurate than expected, namely
second-order. The above proof extends this result to both arbitrary
meshes, and to any finite element pair that satisfies the embedding
properties above.

\begin{figure}
 \begin{center}
   \includegraphics*[width=10cm]{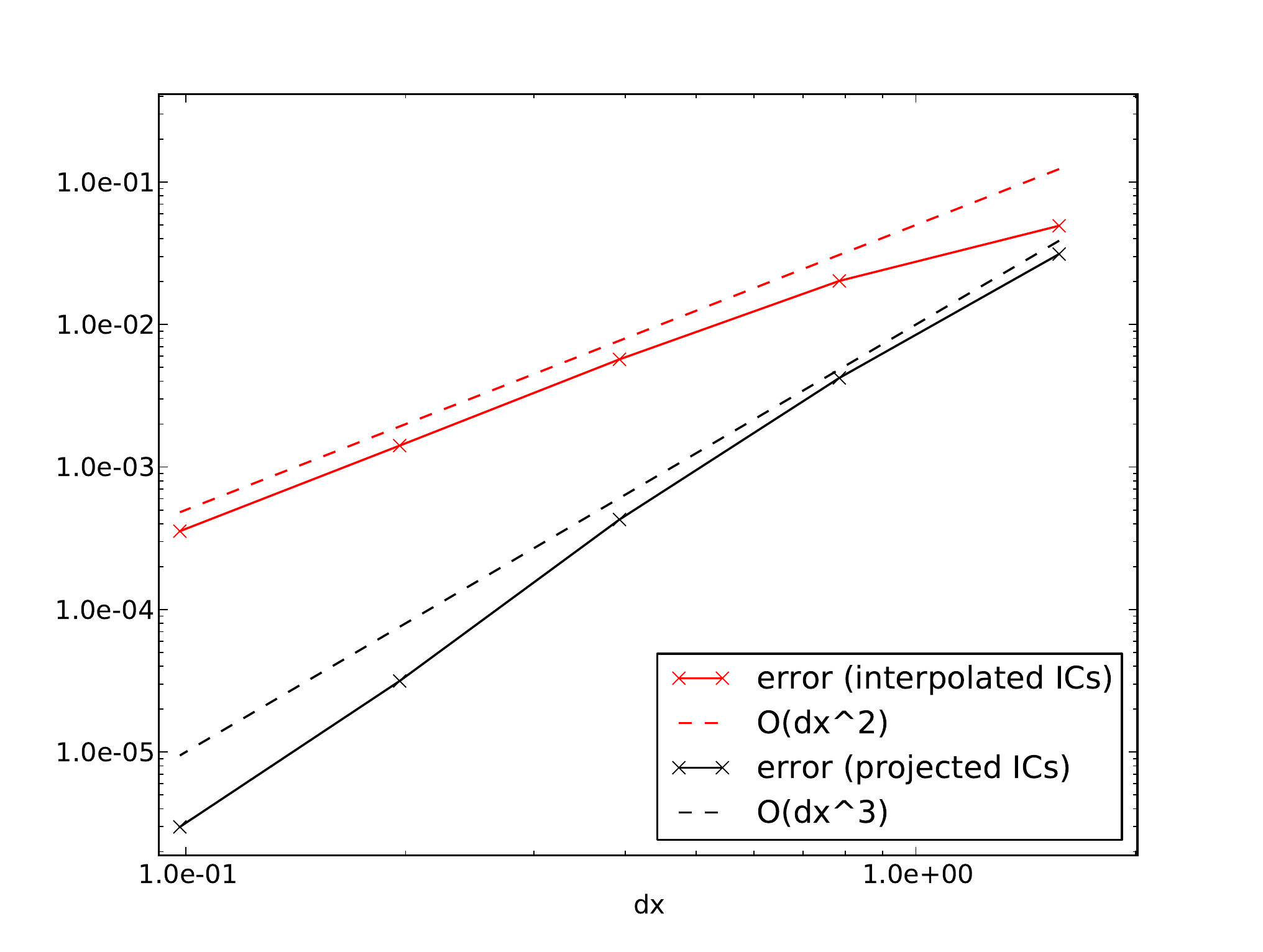}
 \end{center}
 \caption{\label{third-order}Plot showing convergence rates for the
   \pdgp discretisation applied to the linear rotating shallow-water
   equations on an $f$-plane.  The test problem is a single
   propagating sinusoidal wave in periodic boundary conditions, with
   the $L_2$ error in the free surface elevation computed after the
   wave has propagated all the way around the periodic
   domain. Second-order convergence for the free surface elevation is
   obtained when the initial conditions for velocity are obtained by
   collocation at node points; third-order convergence for the free
   surface elevation is obtained when the initial conditions are
   obtained by collocation with a quadratic $P2$ basis for velocity
   and reduced to the \pdg space by $L_2$ projection.}
\end{figure}

A numerical verification of this third-order convergence is shown in
Figure \ref{third-order}. Care must be taken to obtain third-order
convergence: if the initial conditions for the \pdg velocity are
obtained by \pdg collocation, \emph{i.e.} evaluating the analytic
initial condition at the node points and using those values as nodal
basis coefficients, then the truncation error in the initial condition
for the velocity is second-order, and hence second-order accuracy is
the most that can be expected after time-integrating the
equations. However, a third-order accurate velocity initial condition
can be obtained by first constructing a higher-order finite element
approximation to the velocity field by collocation (we used a $P2$
approximation in the calculations in Figure \ref{third-order}), and
then applying the $L_2$ projection to obtain a \pdg velocity
field. This results in third-order convergence of the free surface
elevation over fixed time, since the free surface elevation equation
is the $P2$ finite element approximation to the inertia-gravity wave
equation, as shown above.  To see that this procedure leads to a
third-order accurate velocity field initial condition, first write the
analytic initial condition for the velocity as
\[
\MM{u}(\MM{x},0) = \nabla\phi_0 + \nabla^\perp\psi_0 + \bar{\U}_0.
\]
By standard approximation theory, the $p$th-order collocated finite
element approximation to the initial condition satisfies
$\MM{u}^p=\MM{u}(\MM{x},0)+\mathcal{O}(\Delta x^{p+1})$. The \pdgp
initial condition $\MM{u}^\delta$ satisfies
\[
\int\MM{v}^\delta\cdot\MM{u}^\delta\diff{V}=
\int\MM{v}^\delta\cdot\MM{u}^p\diff{V}
\]
for all \pdg test functions $\MM{v}^\delta$. After subsitution of the
Helmholtz decomposition for $\MM{u}(\MM{x},0)$ and the discrete
Helmholtz decomposition for $\MM{u}^\delta$, this becomes
\begin{eqnarray*}
\int \nabla\alpha^\delta\cdot\nabla\phi^\delta\diff{V} &=&
\int  \nabla\alpha^\delta\cdot\nabla\phi_0\diff{V} +
\mathcal{O}(\Delta x^{p+1}), \\
\int \nabla\alpha^\delta\cdot\nabla\psi^\delta\diff{V} &=&
\int  \nabla\alpha^\delta\cdot\nabla\psi_0\diff{V} +
\mathcal{O}(\Delta x^{p+1}),  \\
\bar{\MM{u}}^\delta & = & \bar{\MM{u}}_0  +
\mathcal{O}(\Delta x^{p+1}), 
\end{eqnarray*}
and the potentials $\phi^\delta$ and $\psi^\delta$ converge to
$\phi_0$ and $\psi_0$ as $\mathcal{O}(\Delta x^3)$ following standard
convergence theory for finite element discretisations of elliptic
problems (see \citet{BrSc04}, for example). Third-order convergence for
the \pdgp discretisation applied to inertia-gravity waves on the
$f$-plane was demonstrated in \citet{CoLaReLe2010} in which various
partly-discontinuous finite element pairs were benchmarked against a
high-order discontinuous Galerkin reference solution. Since the
initial conditions were obtained by $L_2$ projection from the
high-order solution, third-order convergence was observed.

\subsection{Physical inertial oscillation}
\label{phys inertial}
Since the integration is performed over spatially-independent
functions, equation \eqref{f plane inertial} may be written as 
\[
\bar{\W}^\delta\cdot\left(\dd{}{t}\bar{\MM{u}}^\delta +
  f(\bar{\MM{u}}^\delta)^\perp\right)=0,
\]
and since it must hold for all $\bar{\MM{w}}^\delta$, we obtain
\[
\dd{}{t}\bar{\MM{u}}^\delta + f(\bar{\MM{u}}^\delta)^\perp = 0,
\]
which is the usual inertial oscillation equation which has
spatially-independent solutions which rotate with frequency $f$.

\subsection{Spurious inertial oscillations}
\label{spurious inertial}
Equation \eqref{f plane spurious} describes the dynamics of the
spurious velocity component $\hat{\U}^\delta$. If $\hat{\U}^\delta$ is a spurious
velocity (\emph{i.e.} is orthogonal to $\nabla\alpha^\delta$ and
$\nabla^\perp\alpha^\delta$), then so is $(\hat{\U}^\delta)^\perp$ and so
equation \eqref{f plane spurious} does not involve any projection and
hence can be written as
\[
\dd{}{t}\hat{\U}^\delta + f(\hat{\U}^\delta)^\perp = 0.
\]
these solutions also simply rotate with frequency $f$ and hence must
be interpreted as spurious inertial oscillations which do not
propagate as waves.

If we replace the velocity space \pdg with the restricted space \hptwo,
as described in section \ref{helmholtz}, then we obtain the finite
element pair which we call \hptwo-\ptwo. we still have 
equations (\ref{f plane phi}-\ref{f plane inertial}) but without the
spurious inertial oscillations in equation \eqref{f plane spurious},
hence the \hptwo-\ptwo discretisation has no spurious modes.

\subsection{Discrete dispersion relation for inertia-gravity waves}
In this section, we compute the discrete dispersion relation for the
\pdgp discretisation applied to the rotating shallow-water equations
on the $f$-plane for the special case of a structured mesh in a
regular hexagonal domain with edge length $L$ centred on the origin,
with periodic boundary conditions for opposing faces, tiled with
equilateral triangles with edge lengths $\Delta x=L/N$ for some
positive integer $N$, and use this to define a continuous $P2$ finite
element mesh. The discrete dispersion relation is developed by
searching for time-harmonic solutions of \eqref{ig equation}. Assuming
such a time-harmonic solution $\eta^\delta\propto e^{i\omega t}$,
equation \eqref{ig equation} becomes
\begin{equation}
\left(-\omega^2
+f^2\right)
\LL \alpha^\delta, \eta^\delta \RR
+c^2\LL \nabla\alpha^\delta, \nabla\eta^\delta \RR
= 0. \label{ft ig}
\end{equation}
If $\eta^\delta$ is an eigensolution of equation \eqref{ft ig}, then
so is $T_{\MM{z}}\eta^\delta(\MM{x})=\eta^\delta(\MM{x}-\MM{z})$ for
any $\MM{z}$ in the set $\mathbb{V}$ of translations that map vertices
in the mesh to other vertices. Hence, eigenfunctions of equation
\eqref{ft ig} are all eigenfunctions of $T_{\MM{z}}$, \emph{i.e.}
they take the form
\begin{equation}
\label{h form}
\eta^{\delta}(\MM{x})|_{\MM{x}\in\Omega_{\MM{z}}} =
\hat{\eta}^{\delta}(\MM{\xi})e^{i\MM{k}\cdot\MM{z}}, \quad \MM{\xi}\Delta
x + \MM{z} = \MM{x}, \quad \forall\MM{z}\in\mathbb{V},
\end{equation}
where $\Omega_{\MM{z}}$ is the translation of the hexagon formed from
the six equilateral triangles surrounding the vertex at the origin by
$\MM{z}$, $\hat{\eta}^\delta(\MM{\xi})$ is defined on the reference
hexagon $\Omega_e$ with edge length 1 and centred at the origin,
$\MM{\xi}$ is the local coordinate in $\Omega_e$, and
$\MM{k}\in\mathbb{R}^2$ is the wave vector satisfying
$\MM{k}\cdot\MM{z}=2\pi l$ with $l$ an integer. The wave vector
$\MM{k}$ is contained in the first Brillouin zone of the periodic
hexagonal domain which is bounded by the lines
\[
\MM{k}\cdot (\cos(\theta_n),\sin(\theta_n))^T=\frac{2}{\sqrt{3}}\pi,
\quad \theta_n = \left(n+\frac{1}{2}\right)\pi/3, \quad \mbox{for }n=1,2,\ldots,6.
\]
For more details of functions on periodic lattices, see
\citep{Ko2005}, for example.

\begin{figure}
 \begin{center}
   \includegraphics*[width=10cm]{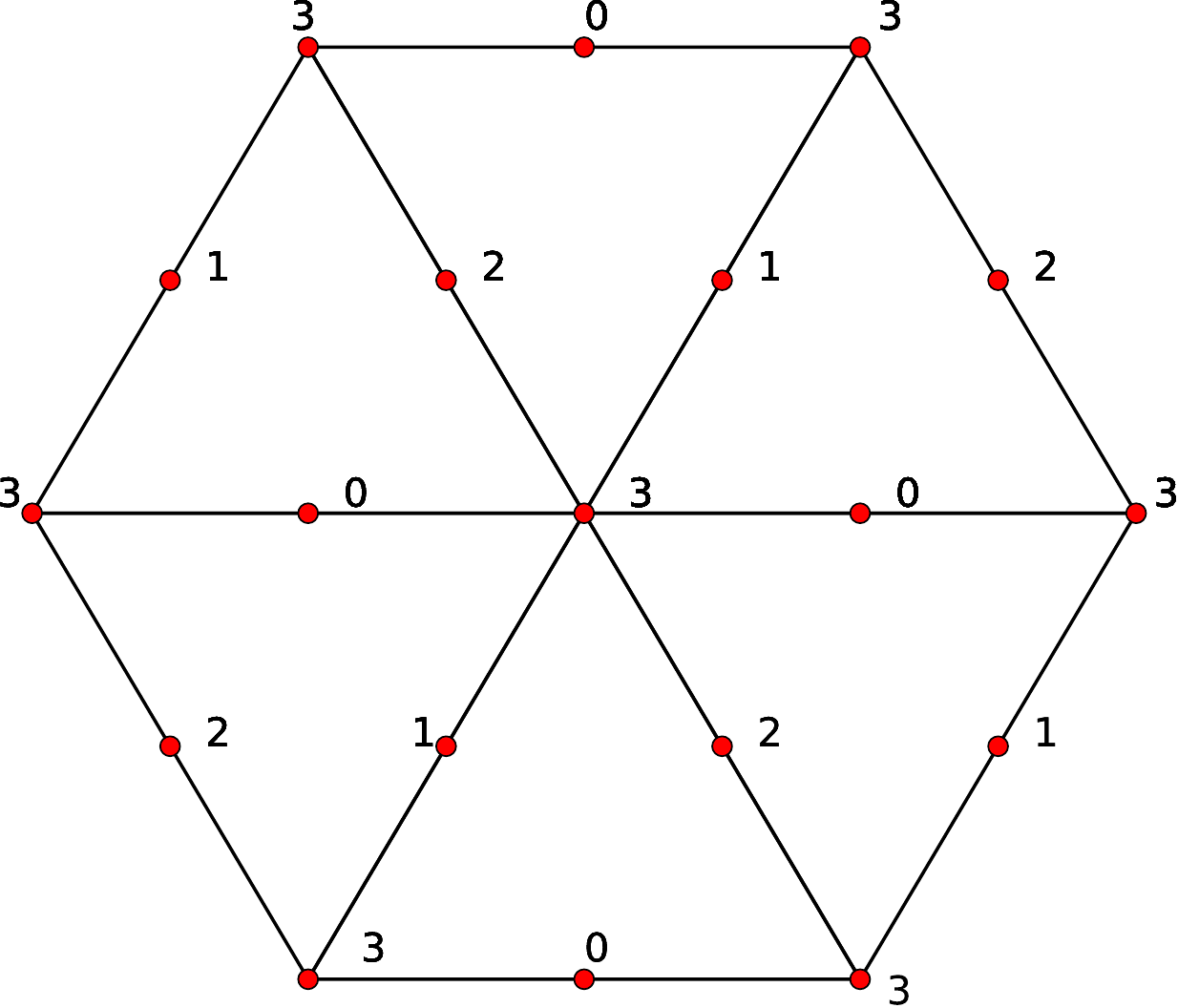}
 \end{center}
 \caption{\label{hexagon} Diagram showing the reference domain
    $\Omega_{e}$ which is used to perform the numerical dispersion
    relation calculations. After considering the boundary conditions
    for $\hat{\eta}^\delta$ which are the consequence of requiring that
    $\eta^\delta$ is continuous, there are four degrees of freedom for
    $\hat{\eta}^\delta$, which we denote $\{\tilde{\eta}_n\}_{n=0}^4$. Each
    node in the diagram is labelled with a number $n$, indicating that
    $\hat{\eta}^\delta = \tilde{\eta}_ne^{i\MM{k}\cdot\MM{\xi}\Delta x}$ at
    that node.}
\end{figure}

Let us now fix an arbitrary wave vector $\MM{k}$ satisfying the
conditions above.  We note that the integral in equation \eqref{ft ig}
can be performed by integrating over all hexagons $\Omega_{\MM{z}}$ and
dividing by three (since each equilateral triangle is covered by three
hexagons).  Given a test function $\alpha^\delta$, equation \eqref{ft
  ig} (multiplied by three) becomes
\begin{eqnarray*}
0 & = & \sum_{\MM{z}\in\mathbb{V}}\int_{\Omega_{\MM{z}}}
\left(-\omega^2
+f^2\right)
\alpha^\delta(\MM{x})\eta^\delta(\MM{x})
+ c^2 \nabla\alpha^\delta(\MM{x})\cdot\nabla\eta^\delta(\MM{x})\diff{V}(\MM{x}) \\
& = & 
\sum_{\MM{z}\in\mathbb{V}}\int_{\Omega_{e}}
\left(\Delta x^2\left(-\omega^2
+f^2\right)
\alpha^\delta(\MM{\xi}\Delta x+\MM{z})\hat{\eta}^\delta(\MM{\xi})
+\nabla_{\MM{\xi}}
\alpha^\delta(\MM{\xi}\Delta x+\MM{z})\cdot\nabla_{\MM{\xi}}
\hat{\eta}^\delta(\MM{\xi})\right)
e^{i\MM{k}\cdot \MM{z}}
\diff{V}(\MM{\xi}), \\
& = & 
\int_{\Omega_{e}}\Delta x^2\left(-\omega^2
+f^2\right)
\hat{\alpha}^\delta(\MM{\xi}) \hat{\eta}^\delta(\MM{\xi})
+ c^2\nabla\hat{\alpha}^\delta(\MM{\xi})\cdot\nabla\hat{\eta}^\delta(\MM{\xi})
\diff{V}(\MM{\xi}), \\
\end{eqnarray*}
where $\hat{\alpha}^\delta$ is defined on $\Omega_e$ with
\[
\hat{\alpha}^\delta(\MM{\xi})=\sum_{\MM{z}\in\mathbb{V}}
\alpha^\delta(\MM{\xi}\Delta x-\MM{z})
e^{i\MM{k}\cdot\MM{z}}.
\]
We have now written the dispersion relation in such a way that all the
computations can be done over one single reference hexagon
$\Omega_{e}$.  The boundary conditions for $\hat{\eta}^\delta$ on the
reference hexagon can be computed from the condition that
$\hat{\eta}^\delta$ is continuous at the boundaries, meaning that on each
edge of the hexagon $\Omega_{e}$, denoted $\partial\Omega_{e,n}$ (with
$1\leq n\leq6$),
\[
\hat{\eta}^\delta(\MM{\xi}) =
e^{i\Delta x\MM{k}\cdot\Delta\MM{\xi}}\hat{\eta}^{\delta}(\MM{\xi}+\Delta\MM{\xi}),
\]
where $\Delta\MM{\xi}$ is the vector from $\partial\Omega_{e,n}$ to
the opposing face. Figure \ref{hexagon} illustrates the consequences
of this for the basis coefficients of $\hat{\eta}^{\delta}$ when a
nodal basis\footnote{A nodal basis is a basis in which each basis
  function has unit value at one of the node points, \emph{e.g.} the
  vertices and edge midpoints in the case of the continuous quadratic
  mesh, and vanishes on all other node points.} is used.

We can similarly use continuity of $\alpha^\delta$ to obtain boundary
conditions for $\hat{\alpha}^\delta(\MM{\xi})$ on
$\partial\Omega_e$. On the boundary $\partial\Omega_{e,n}$, 
\begin{eqnarray*}
  \hat{\alpha}^\delta(\MM{\xi})
  & = & \sum_{\MM{z}\in\mathbb{V}}
  \alpha^\delta(\MM{\xi}\Delta x-\MM{z})
  e^{i\MM{k}\cdot\MM{z}}, \\
  & = & \sum_{\MM{z}\in\mathbb{V}}
  \alpha^\delta((\MM{\xi}+\Delta\MM{\xi})\Delta x-(\MM{z}+\Delta x\Delta\MM{\xi}))
  e^{i\MM{k}\cdot\MM{z}}, \\
  & = & \sum_{\MM{z}\in\mathbb{V}}
  \alpha^\delta((\MM{\xi}-\Delta\MM{\xi})\Delta x-\MM{z})
  e^{i\MM{k}\cdot(\MM{z}-\Delta x\Delta\MM{\xi})}, \\
  & = & e^{-i\Delta x\MM{k}\cdot\MM{\xi}}
\hat{\alpha}^\delta(\MM{\xi}+\Delta\MM{\xi}).
\end{eqnarray*}
This means that $\hat{\alpha}^\delta$ has boundary conditions which
are the complex conjugate of the boundary conditions for
$\hat{\eta}^\delta$.

We adopt a nodal basis for functions inside $\Omega_e$. There are 
19 P2 nodes on $\Omega_e$ (see Figure \ref{hexagon}), and so we write
\[
\hat{\eta}^\delta = \sum_{n=1}^{19}\hat{\eta}_nN_n(\MM{\xi}),
\]
where $N_n(\MM{\xi})$, $(n=1,\ldots,19)$, are the nodal basis
functions for P2 functions inside $\Omega_e$, and $\eta_n$
$(n=1,\ldots,19)$ are the nodal basis coefficients.  The boundary
conditions for $\hat{\eta}^\delta$ described above can be expressed
\emph{via} a matrix $S$ (which is a function of $\MM{k}\Delta x$ due
to the dependence of the boundary conditions for $\hat{\eta}$ and
$\hat{\alpha}$ on $\MM{k}$), so that
\[
\MM{\hat{\eta}} = S\MM{\tilde{\eta}}, \quad \MM{\hat{\alpha}} =
S^*\MM{\tilde{\alpha}},
\]
where $\MM{\hat{\eta}}$ and $\MM{\hat{\alpha}}$ are the vectors of the
basis coefficients of $\eta^{\delta}$ and $\alpha^{\delta}$
respectively, and $\MM{\tilde{\eta}}$ and $\MM{\tilde{\alpha}}$ are the
corresponding vectors of the independent degrees of freedom. 

After substituting, the wave equation becomes
\begin{eqnarray*}
0 & = & \left(-\omega^2
+f^2\right)
\LL \alpha^\delta, \eta^\delta \RR
+g\LL \nabla\alpha^\delta, \nabla\eta^\delta \RR
= 0 \\
& = & 
\Delta x^2\hat{\MM{\alpha}}^T
\left(
\left(-\omega^2
+f^2\right)M_e
+\frac{g L_e}{\Delta x^2}\right)\hat{\MM{\eta}} \\
& = & \Delta x^2\tilde{\MM{\alpha}}^TS^{\dagger}
\left(
\left(-\omega^2
+f^2\right)M_e
+g \frac{L_e}{\Delta x^2}\right)S\tilde{\MM{\eta}},
\end{eqnarray*}
where $M_e$ is the local mass matrix
\[
M_{e,ij} = \int_{\Omega_e}N_i(\MM{\xi})N_j(\MM{\xi})\diff{V}(\MM{\xi}),
\]
and $L_e$ is the Laplacian matrix
\[
L_{e,ij}=\int_{\Omega_e}\nabla N_i(\MM{\xi})\cdot 
\nabla N_j(\MM{\xi})\diff{V}(\MM{\xi}),
\]
and $\dagger$ indicates the Hermitian conjugate of a matrix. Since
$\tilde{\MM{\alpha}}$ is arbitrary, we seek non-trivial solutions of
\[
S^{\dagger}
\left(
\Delta x^2\left(-\omega^2
+
f^2\right)M_e
+g L_e\right)S\tilde{\MM{\eta}} = \MM{0},
\]
and we obtain the dispersion relation
\begin{equation}
\label{ig dispersion relation}
\left|S^{\dagger}
\left(\Delta x^2
\left(-\omega^2
+f^2\right)M_e
+g L_e\right)S\right|=0,
\end{equation}
which must be solved for $\omega$ given $\MM{k}$ (the $\MM{k}$
dependence is in $S$ as described above). This equation is the
determinant of a $4\times 4$ matrix with entries that are linear in
$\lambda=\Delta x^2(\omega^2-f^2)$, so it is quartic polynomial in
$\lambda$.

After lengthy calculation using SymPy \citep{SymPy}, the
following matrices are obtained:
\[
S^\dagger M_eS = \begin{pmatrix}
A & B \\
B^T & C \\
\end{pmatrix}, \qquad
S^\dagger L_e S = \begin{pmatrix}
D & E \\
E^T & F \\
\end{pmatrix},
\]
where
\begin{eqnarray*}
  A&=&\left(\begin{smallmatrix}\frac{4}{15} \sqrt{3} & \frac{2}{15}
      \sqrt{3} \operatorname{cos}\left(- \frac{1}{4} k + \frac{1}{4} l
        \sqrt{3}\right)\\\frac{2}{15} \sqrt{3} \operatorname{cos}\left(-
        \frac{1}{4} k + \frac{1}{4} l \sqrt{3}\right) & \frac{4}{15}
      \sqrt{3}\end{smallmatrix}\right), \\
  B &=& \left(\begin{smallmatrix}\frac{2}{15} \sqrt{3}
      \operatorname{cos}\left(\frac{1}{4} k + \frac{1}{4} l
        \sqrt{3}\right) & \frac{2}{15} \sqrt{3}
      \operatorname{cos}\left(\frac{1}{2} k\right)\\- \frac{1}{30}
      \sqrt{3} \operatorname{cos}\left(\frac{1}{2} l \sqrt{3}\right) & -
      \frac{1}{30} \sqrt{3} \operatorname{cos}\left(\frac{3}{4} k -
        \frac{1}{4} l \sqrt{3}\right)\end{smallmatrix}\right), \\
  C & = & \left(\begin{smallmatrix}\frac{4}{15} \sqrt{3} & -
      \frac{1}{30} \sqrt{3} \operatorname{cos}\left(\frac{3}{4} k +
        \frac{1}{4} l \sqrt{3}\right)\\- \frac{1}{30} \sqrt{3}
      \operatorname{cos}\left(\frac{3}{4} k + \frac{1}{4} l
        \sqrt{3}\right) & - \frac{1}{60} \sqrt{3}
      \operatorname{cos}\left(k\right) - \frac{1}{60} \sqrt{3}
      \operatorname{cos}\left(\frac{1}{2} k + \frac{1}{2} l
        \sqrt{3}\right) - \frac{1}{60} \sqrt{3}
      \operatorname{cos}\left(- \frac{1}{2} k + \frac{1}{2} l
        \sqrt{3}\right) + \frac{3}{20} \sqrt{3}\end{smallmatrix}\right),
  \\
  D & = &   \left(\begin{smallmatrix}8 \sqrt{3} & - \frac{8}{3} \sqrt{3}
      \operatorname{cos}\left(- \frac{1}{4} k + \frac{1}{4} l
        \sqrt{3}\right)\\- \frac{8}{3} \sqrt{3}
      \operatorname{cos}\left(- \frac{1}{4} k + \frac{1}{4} l
        \sqrt{3}\right) & 8 \sqrt{3}\end{smallmatrix}\right), \\
  E & = & \left(\begin{smallmatrix}- \frac{8}{3} \sqrt{3}
      \operatorname{cos}\left(\frac{1}{4} k + \frac{1}{4} l
        \sqrt{3}\right) & - \frac{8}{3} \sqrt{3}
      \operatorname{cos}\left(\frac{1}{2} k\right)\\- \frac{8}{3}
      \sqrt{3} \operatorname{cos}\left(\frac{1}{2} k\right) & -
      \frac{8}{3} \sqrt{3} \operatorname{cos}\left(\frac{1}{4} k +
        \frac{1}{4} l \sqrt{3}\right)\end{smallmatrix}\right),
  \quad\mbox{and} \\
F & = & \left(\begin{smallmatrix}8 \sqrt{3} & - \frac{8}{3} \sqrt{3} \operatorname{cos}\left(- \frac{1}{4} k + \frac{1}{4} l \sqrt{3}\right)\\- \frac{8}{3} \sqrt{3} \operatorname{cos}\left(- \frac{1}{4} k + \frac{1}{4} l \sqrt{3}\right) & \frac{2}{3} \sqrt{3} \operatorname{cos}\left(- \frac{1}{2} k + \frac{1}{2} l \sqrt{3}\right) + \frac{2}{3} \sqrt{3} \operatorname{cos}\left(k\right) + \frac{2}{3} \sqrt{3} \operatorname{cos}\left(\frac{1}{2} k + \frac{1}{2} l \sqrt{3}\right) + 6 \sqrt{3}\end{smallmatrix}\right),
\end{eqnarray*}
having written $\MM{k}=(k,l)$.

\begin{figure}
 \begin{center}
   \includegraphics*[width=10cm]{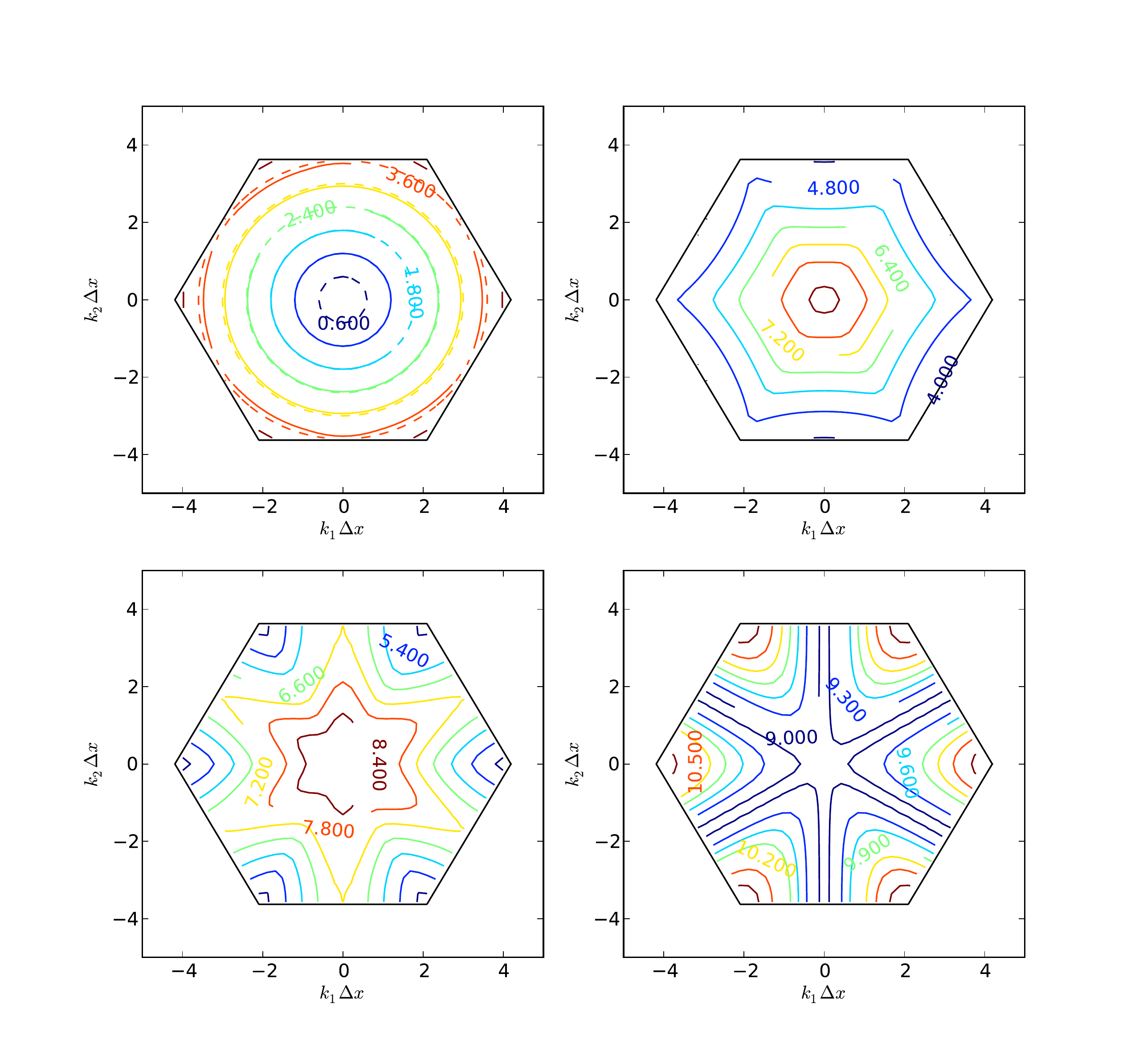}
 \end{center}
 \caption{\label{poincare_spec} Plots showing contours of $\Delta
   x(\omega^2-f^2)$ in the $\MM{k}\Delta x$ plane for each of the
   four branches of the numerical dispersion relation for the \pdgp
   finite element scheme applied to the linear rotating shallow water
   equations on the $f$-plane. The lowest branch is shown top-left,
   with contours of the exact dispersion relation superimposed using
   dashed lines. This lowest branch is very accurate, and the
   contours are very circular, meaning that the wave propagation is
   almost independent of the direction of mesh alignment.  The other
   plots show the higher branches which represent the second, third
   and fourth Brillouin zones in the $\MM{k}\Delta x$ plane mapped in
   to the first Brillouin zone. For example, one can cross from the
   lowest branch into the branch in the top-right branch by going
   through the hexagon which bounds the region, emerging from the
   opposite edge in the hexagon in the top-right plot, moving in the
   opposite direction. It can be seen that all four branches
   represent physical modes from different regions of physical
   $\MM{k}$-space which can be resolved on the grid.}
\end{figure}

The resulting quartic equation for $\lambda=\Delta x^2(\omega^2-f^2)$
obtained from evaluating the determinant \eqref{ig dispersion
  relation} is a very complicated expression that would take up
several pages. Hence, solutions to the dispersion relation equation
\eqref{ig dispersion relation} were obtained by numerically evaluating
the matrix $(S^\dagger M_e S)^{-1} S^\dagger L_e S$ for various values
of $\MM{k}$, and using the Scientific Python {\ttfamily linalg.eig}
routine, which were then sorted in numerical order. Since the equation
for $\lambda=\Delta x^2(\omega^2-f^2)$ is quartic, this leads to four
branches of the dispersion relation (this is typical for P2 schemes in
two dimensions), which correspond to the fundamental
$\exp(i\MM{k}\cdot\MM{x})$ modes with $\MM{k}$ inside the first
Brillouin zone, plus higher wave number solutions obtained from the
second, third and fourth Brillouin zones which have the same
translation property at the triangle vertices but result in different
values at the edge centres. The plots of the four branches are given
in Figure \ref{poincare_spec}. It is immediately visible that the
lowest eigenvalues are very isotropic, as might be expected from the
fact that the dispersion relation is in fact third-order rather than
second-order, as described in section \ref{inertia gravity
  waves}. This means that resolved gravity waves of a particular wave
number have a propagation speed which is largely independent of the
direction of alignment of the mesh (this is a property which is
considered important and was one of the contributing factors towards
designing the hexagonal C-grid as an alternative to the triangular
C-grid). It can also be seen that the dispersion relation is
monotonically-increasing with $|\MM{k}|$ with some small jumps when
moving between branches (see \citet{CoHaPaRe2009} for the equivalent
one-dimensional plot); there are no spurious inertia-gravity modes.

\section{Discrete wave propagation on the $\beta$-plane}
\label{beta plane}
In this section, we consider the quasi-geostrophic scaling on the
$\beta$-plane, following the approach of \citet{LePo2008,Th08} in which the
quasi-geostrophic approximation is applied to the
spatially-discretised equations.

In the $\beta$-plane case, $f=f_0+\beta y$, and after substitution 
of the orthogonal decomposition for the solution variables and test
functions into equations (\ref{u eqn}-\ref{h eqn})
 we obtain
\begin{eqnarray*}
\dd{}{t}\LL \nabla\alpha^\delta, \nabla \phi^\delta \RR
- f_0\LL \nabla\alpha^\delta, \nabla\psi^\delta \RR 
- \LL \beta y\nabla\alpha^\delta, \bar{\U}^\delta+ (\hat{\U}^\delta)^\perp + 
\nabla\psi^\delta + \nabla^\perp \phi^\delta
\RR
+ c^2\LL \nabla\alpha^\delta, \nabla\eta^\delta \RR
& = & 0  \\
\dd{}{t}\LL\nabla\alpha^\delta, \nabla\psi^\delta\RR 
+ f_0\LL \nabla\alpha^\delta,\nabla\phi^\delta \RR 
+ \LL \beta y\nabla\alpha^\delta,-\bar{\U}^\delta-\hat{\U}^\delta+
\nabla^\perp\phi^\delta + \nabla\psi^\delta \RR 
& = & 0 \\
\dd{}{t}\LL \alpha^\delta, \eta^\delta \RR
 - \LL\nabla\alpha^\delta,\nabla\phi^\delta\RR & = & 0, \\
\dd{}{t} \LL \bar{\W}^\delta,\bar{\U}^\delta \RR +
f_0 \LL \bar{\W}^\delta, (\bar{\U}^\delta)^\perp \RR
+ \LL \bar{\W}^\delta\beta y ,(\bar{\U}^\delta)^\perp+(\hat{\U}^\delta)^\perp +
\nabla^\perp\phi^\delta - \nabla\psi^\delta\RR & = & 0, \\
\dd{}{t} \LL \hat{\W}^\delta,\hat{\U}^\delta \RR +
f_0 \LL \hat{\W}^\delta, (\hat{\U}^\delta)^\perp \RR
+ \LL \hat{\W}^\delta\beta y ,(\bar{\U}^\delta)^\perp+ (\hat{\U}^\delta)^\perp+ 
\nabla^\perp\phi^\delta - \nabla\psi^\delta\RR & = & 0.
\end{eqnarray*}
At leading order in Rossby number in the quasi-geostrophic scaling, we
obtain the geostrophic balance:
\begin{eqnarray*}
- f_0\LL \nabla\alpha^\delta, \nabla\psi^\delta_g \RR 
+ c^2\LL \nabla\alpha^\delta, \nabla\eta^\delta_g \RR
& = & 0, 
\\
f_0\LL \nabla\alpha^\delta,\nabla\phi^\delta_g \RR 
& = & 0, \\
 - \LL\nabla\alpha^\delta,\nabla\phi^\delta_g\RR & = & 0, \\
f_0\LL \hat{\W}^\delta,(\hat{\U}^\delta)^\perp_g \RR & = & 0,
\end{eqnarray*}
which we have already analysed in Section \ref{geostrophic balance},
and so we know that it implies that
\begin{equation}
\label{balance}
\hat{\U}^\delta_g = 0, \quad \phi^\delta_g = 0, \quad \psi^\delta_g = \frac{c^2}{f}\eta^\delta_g.
\end{equation}
At the next order we obtain
\begin{eqnarray}
\label{phi eqn qg}
\dd{}{t}\LL \nabla\alpha^\delta, \nabla \phi^\delta_{ag} \RR
- f_0\LL \nabla\alpha^\delta, \nabla\psi^\delta_{ag} \RR 
- \LL \beta y\nabla\alpha^\delta, \nabla\psi^\delta_g\RR
+ gH\LL \nabla\alpha^\delta, \nabla\eta^\delta_{ag} \RR
& = & 0  \\
\dd{}{t}\LL\nabla\alpha^\delta, \nabla\psi^\delta_g\RR 
+ f_0\LL \nabla\alpha^\delta,\nabla\phi^\delta_{ag} \RR 
+ \LL \beta y\nabla\alpha^\delta,
\nabla^\perp \psi^\delta_g \RR 
& = & 0 \label{psi eqn qg}
\\
\label{eta eqn qg}
\dd{}{t}\LL \alpha^\delta, \eta^\delta_{g} \RR
 - \LL\nabla\alpha^\delta,\nabla\phi^\delta_{ag} \RR & = & 0, \\
f_0 \LL \bar{\W}^\delta, (\bar{\U}^\delta)^\perp_{ag} \RR
+ \LL \bar{\W}^\delta\beta y , - 
\nabla\psi^\delta_g \RR & = & 0 ,\\
f_0 \LL \hat{\W}^\delta, (\hat{\U}^\delta)^\perp_{ag} \RR
+ \LL \hat{\W}^\delta\beta y , - 
\nabla\psi^\delta_g \RR & = & 0 
\label{spurious qg} 
\end{eqnarray}
Notice that the spurious velocity modes do not appear at this order in
the physical mode equations (\ref{phi eqn qg}-\ref{eta eqn qg}), and
that equation \eqref{spurious qg} states that the ageostrophic
spurious velocity modes are slaved to the geostrophic
streamfunction. Substituting equations \eqref{balance} and 
\eqref{eta eqn qg} into \eqref{psi
  eqn qg} gives
\begin{equation}
\label{rossby 1}
\dd{}{t}\left(
\LL\nabla\alpha^\delta, \nabla\psi^\delta_g\RR 
+ \frac{f_0^2}{gH}\LL \alpha^\delta,\psi^\delta_g \RR 
\right)
+ \LL \beta y\nabla\alpha^\delta,
\nabla^\perp \psi^\delta_g \RR 
 =  0.
\end{equation}
The second term in equation \eqref{rossby 1} may be written as
\[
\LL \beta y\nabla\alpha^\delta,
\nabla^\perp \psi^\delta_g \RR 
= \LL \nabla(\beta y\alpha^\delta)-\beta \alpha^\delta(0,1),
\nabla^\perp \psi^\delta_g \RR
= -\beta \LL \alpha^\delta,\pp{}{x}\psi^\delta_g\RR,
\]
and we obtain the usual continuous finite element approximation
to the Rossby wave equation using $P2$ elements
\begin{equation}
\dd{}{t}\left(
\LL\nabla\alpha^\delta, \nabla\psi^\delta_g\RR 
+ \frac{f_0^2}{gH}\LL \alpha^\delta,\psi^\delta_g \RR 
\right)
- \beta \LL \alpha^\delta,
\pp{}{x}\psi^\delta_g \RR 
 =  0.
\label{discrete rossby}
\end{equation}
Since $P2$ elements are used, the approximation to the Rossby wave
equation is third-order accurate, rather than the second-order
accuracy one would expect with \pdg for velocity. The equivalent
property for $P0$-$P1$ was shown in \citet{LePo2008}, namely that the
Rossby wave dispersion relation is second-order. The above proof
extends this result to arbitrary meshes and to any finite element pair
which satisfies the embedding properties.

We again expect that the phase velocity is more independent of mesh
orientation than other second-order methods. Since the streamfunction
$\psi^\delta$ and the height variable $\eta^\delta$ are both from the
\ptwo space and hence have the same numbers of degrees of freedom,
there are exactly twice as many inertia-gravity wave modes as Rossby
wave modes. We also note that if the reduced space \hptwo-\ptwo is
used instead of \pdgp we obtain the same equations but with vanishing
spurious inertial modes.

\subsection{Discrete dispersion relation for Rossby waves}

Starting from equation \eqref{discrete rossby}, and following the
method described above for the obtaining the inertia-gravity wave
dispersion relation on the equilateral grid, we obtain the numerical
dispersion relation
\begin{equation}
\label{rossby dispersion relation}
\left| S^\dagger
\left(i\omega\left(\frac{L_e}{\Delta x^2} 
+ \frac{1}{L_R^2}M_e\right) - 
\beta \frac{D_e}{\Delta x}\right)S\right|=0,
\end{equation}
where $D_e$ is the local derivative matrix
\[
D_{e,ij} = \int_{\Omega_e}N_i(\MM{\xi})\hat{\MM{f}}.\nabla
N_j(\MM{\xi})\diff{V}(\MM{\xi}),
\]
and where $\hat{\MM{f}}$ is the unit vector pointing in the direction
of increasing $f$ on the $\beta$-plane. We shall investigate the
variation in the dispersion relation with the alignment of the
triangular grid, and hence it is convenient to write
\[
D_{e,ij} = \hat{f}_1D^1_{e,ij} + \hat{f}_2D^2_{e,ij},
\]
where
\[
D^1_{e,ij} = \int_{\Omega_e}N_i(\MM{\xi})
\pp{N_j}{\xi_1}(\MM{\xi})\diff{V}(\MM{\xi}),\quad
D^2_{e,ij} = \int_{\Omega_e}N_i(\MM{\xi})
\pp{N_j}{\xi_2}(\MM{\xi})\diff{V}(\MM{\xi}).
\]
After further algebraic manipulation with SymPy, we obtain
\[
S^\dagger D_e^1S = 
\begin{pmatrix} 
P_1 & Q_1 \\ 
Q^T_1 & R_1 \\
\end{pmatrix}, \qquad 
S^\dagger L_e S = 
\begin{pmatrix}
P_2 & Q_2 \\
Q^T_2 & R_2 \\
\end{pmatrix},
\]
where
\begin{eqnarray*}
P_1 & = & \left(\begin{smallmatrix}0 & - \frac{2}{5} \sqrt{3}
    \operatorname{sin}\left(- \frac{1}{4} k + \frac{1}{4} l
      \sqrt{3}\right)\\- \frac{2}{5} \sqrt{3}
    \operatorname{sin}\left(- \frac{1}{4} k + \frac{1}{4} l
      \sqrt{3}\right) & 0\end{smallmatrix}\right), \\
Q_1 & = & \left(\begin{smallmatrix}\frac{2}{5} \sqrt{3}
    \operatorname{sin}\left(\frac{1}{4} k + \frac{1}{4} l
      \sqrt{3}\right) & \frac{4}{5} \sqrt{3}
    \operatorname{sin}\left(\frac{1}{2} k\right)\\\frac{3}{5} \sqrt{3}
    \operatorname{sin}\left(\frac{1}{2} k\right) & - \frac{1}{10}
    \sqrt{3} \operatorname{sin}\left(\frac{3}{4} k - \frac{1}{4} l
      \sqrt{3}\right) + \frac{3}{10} \sqrt{3}
    \operatorname{sin}\left(\frac{1}{4} k + \frac{1}{4} l
      \sqrt{3}\right)\end{smallmatrix}\right), \\
R_1 & = & \left(\begin{smallmatrix}0 & - \frac{3}{10} \sqrt{3}
    \operatorname{sin}\left(- \frac{1}{4} k + \frac{1}{4} l
      \sqrt{3}\right) - \frac{1}{10} \sqrt{3}
    \operatorname{sin}\left(\frac{3}{4} k + \frac{1}{4} l
      \sqrt{3}\right)\\- \frac{3}{10} \sqrt{3}
    \operatorname{sin}\left(- \frac{1}{4} k + \frac{1}{4} l
      \sqrt{3}\right) - \frac{1}{10} \sqrt{3}
    \operatorname{sin}\left(\frac{3}{4} k + \frac{1}{4} l
      \sqrt{3}\right) & - \frac{1}{5} \sqrt{3}
    \operatorname{sin}\left(k\right) - \frac{1}{10} \sqrt{3}
    \operatorname{sin}\left(\frac{1}{2} k + \frac{1}{2} l
      \sqrt{3}\right) - \frac{1}{10} \sqrt{3}
    \operatorname{sin}\left(\frac{1}{2} k - \frac{1}{2} l
      \sqrt{3}\right) \end{smallmatrix}\right), \\
P_2 & = & \left(\begin{smallmatrix}0 & \frac{6}{5}
    \operatorname{sin}\left(- \frac{1}{4} k + \frac{1}{4} l
      \sqrt{3}\right)\\\frac{6}{5} \operatorname{sin}\left(-
      \frac{1}{4} k + \frac{1}{4} l \sqrt{3}\right) &
    0\end{smallmatrix}\right), \\
Q_2 & = & \left(\begin{smallmatrix}\frac{6}{5}
    \operatorname{sin}\left(\frac{1}{4} k + \frac{1}{4} l
      \sqrt{3}\right) & 0\\- \frac{1}{5}
    \operatorname{sin}\left(\frac{1}{2} l \sqrt{3}\right) & -
    \frac{1}{10} \operatorname{sin}\left(- \frac{3}{4} k + \frac{1}{4}
      l \sqrt{3}\right) + \frac{9}{10}
    \operatorname{sin}\left(\frac{1}{4} k + \frac{1}{4} l
      \sqrt{3}\right)\end{smallmatrix}\right), \quad \mbox{and} \\
R_2 & = & \left(\begin{smallmatrix}0 & - \frac{1}{10}
    \operatorname{sin}\left(\frac{3}{4} k + \frac{1}{4} l
      \sqrt{3}\right) + \frac{9}{10} \operatorname{sin}\left(-
      \frac{1}{4} k + \frac{1}{4} l \sqrt{3}\right)\\- \frac{1}{10}
    \operatorname{sin}\left(\frac{3}{4} k + \frac{1}{4} l
      \sqrt{3}\right) + \frac{9}{10} \operatorname{sin}\left(-
      \frac{1}{4} k + \frac{1}{4} l \sqrt{3}\right) & - \frac{3}{10}
    \operatorname{sin}\left(\frac{1}{2} k + \frac{1}{2} l
      \sqrt{3}\right) - \frac{3}{20} \operatorname{sin}\left(-
      \frac{1}{2} k + \frac{1}{2} l \sqrt{3}\right) + \frac{3}{20}
    \operatorname{sin}\left(\frac{1}{2} k - \frac{1}{2} l
      \sqrt{3}\right)\end{smallmatrix}\right).
\end{eqnarray*}
The eigenvalues can then be obtained using the method used for the
inertia-gravity waves \emph{i.e.} by finding the eigenvalues of the
matrix for various $\MM{k}\Delta x$ and plotting contours in $\MM{k}$
space. There is an extra difficulty in the Rossby case, because the
numerical algorithm for obtaining eigenvalues of the $4\times 4$
matrix does not preserve the order of the branches when $\MM{k}\Delta
x$ is varied. It is not possible to distinguish the branches by
sorting the eigenvalues in numerical order for each $\MM{k}$ because
the branches have values which cross. However, the branches can be
distinguished by examining the corresponding eigenvectors. If we
interpolate the continuous Fourier modes to the reference hexagon, we
obtain four types of solution (after normalisation) for
$\tilde{\MM{\eta}}$, namely
\[
\begin{pmatrix}
\frac{1}{4} \\
\frac{1}{4} \\
\frac{1}{4} \\
\frac{1}{4} \\
\end{pmatrix}, \quad
\begin{pmatrix}
\frac{-1}{4} \\
\frac{-1}{4} \\
\frac{1}{4} \\
\frac{1}{4} \\
\end{pmatrix}, \quad
\begin{pmatrix}
\frac{1}{4} \\
\frac{-1}{4} \\
\frac{-1}{4} \\
\frac{1}{4} \\
\end{pmatrix}, \quad
\begin{pmatrix}
\frac{-1}{4} \\
\frac{1}{4} \\
\frac{-1}{4} \\
\frac{1}{4} \\
\end{pmatrix},
\]
where the fundamental modes take the form of the vector on the left,
and higher modes arise from the other three vectors. Hence, we
identified the various branches by inspecting the eigenvectors and
associating them with the branch which has the same sign pattern as
the vectors above.

\begin{figure}
\begin{center}
\includegraphics*[width=10cm]{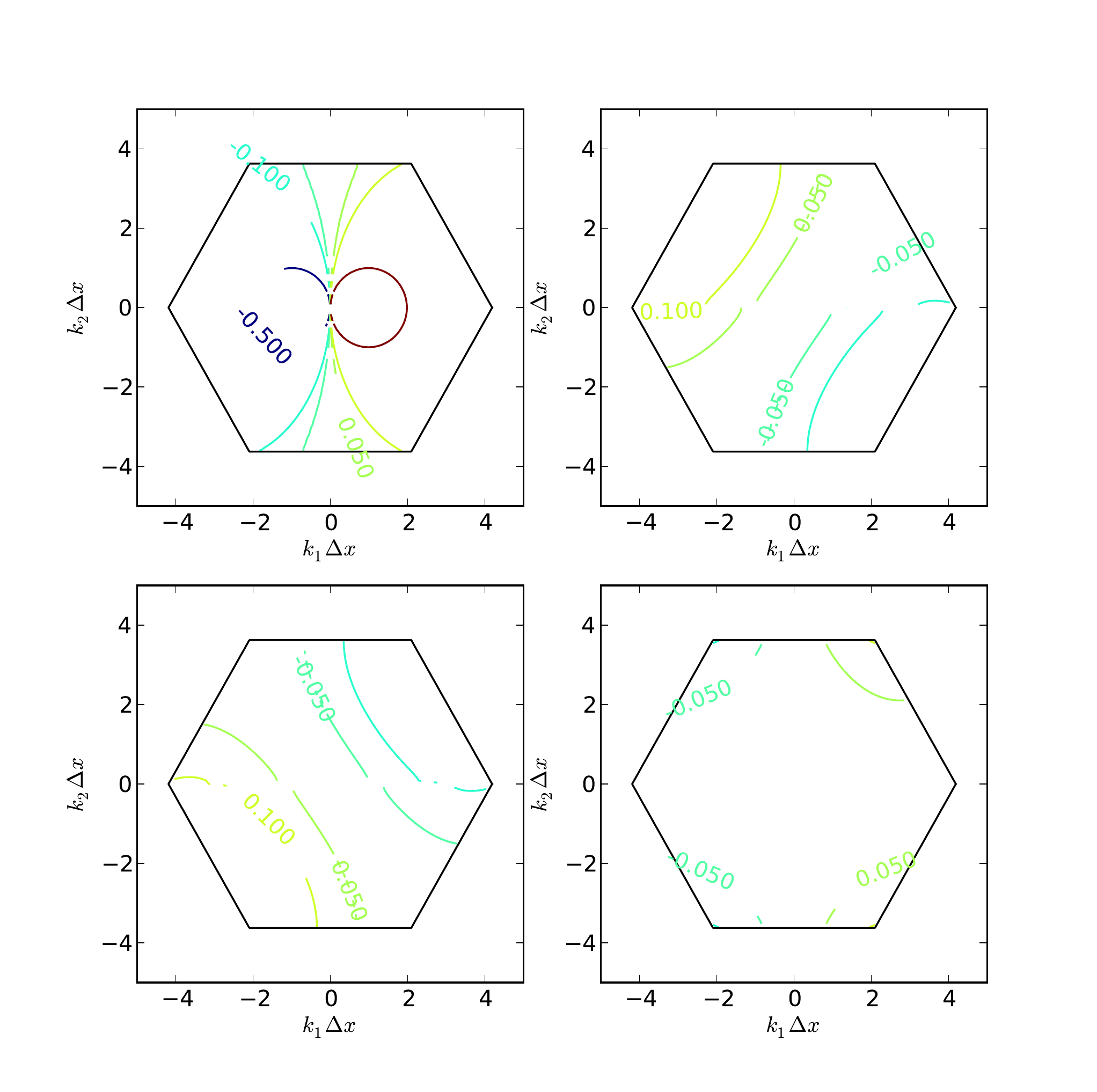}
\end{center}
\caption{\label{rossby_horiz_all_modes} Contour plots showing
  $\omega\times 10^6$ obtained from the solutions of equation
  \eqref{rossby dispersion relation}, with parameters $f_0 =
  1.0\times10^{-4}$, $\beta = 1.0\times 10^{-12}$, $\Delta x =
  1.0\times 10^5$ and $c^2=1.0\times 10^5$ (these parameters are the
  same as those used in \citet{Th08}). $f$ increases in the
  $y$-direction relative to the mesh. The lowest branch of the
  dispersion relation is shown top-left. The other branches are
  aliased higher values of $\MM{k}\Delta x$.}
\end{figure}

\begin{figure}
\begin{center}
\includegraphics*[width=10cm]{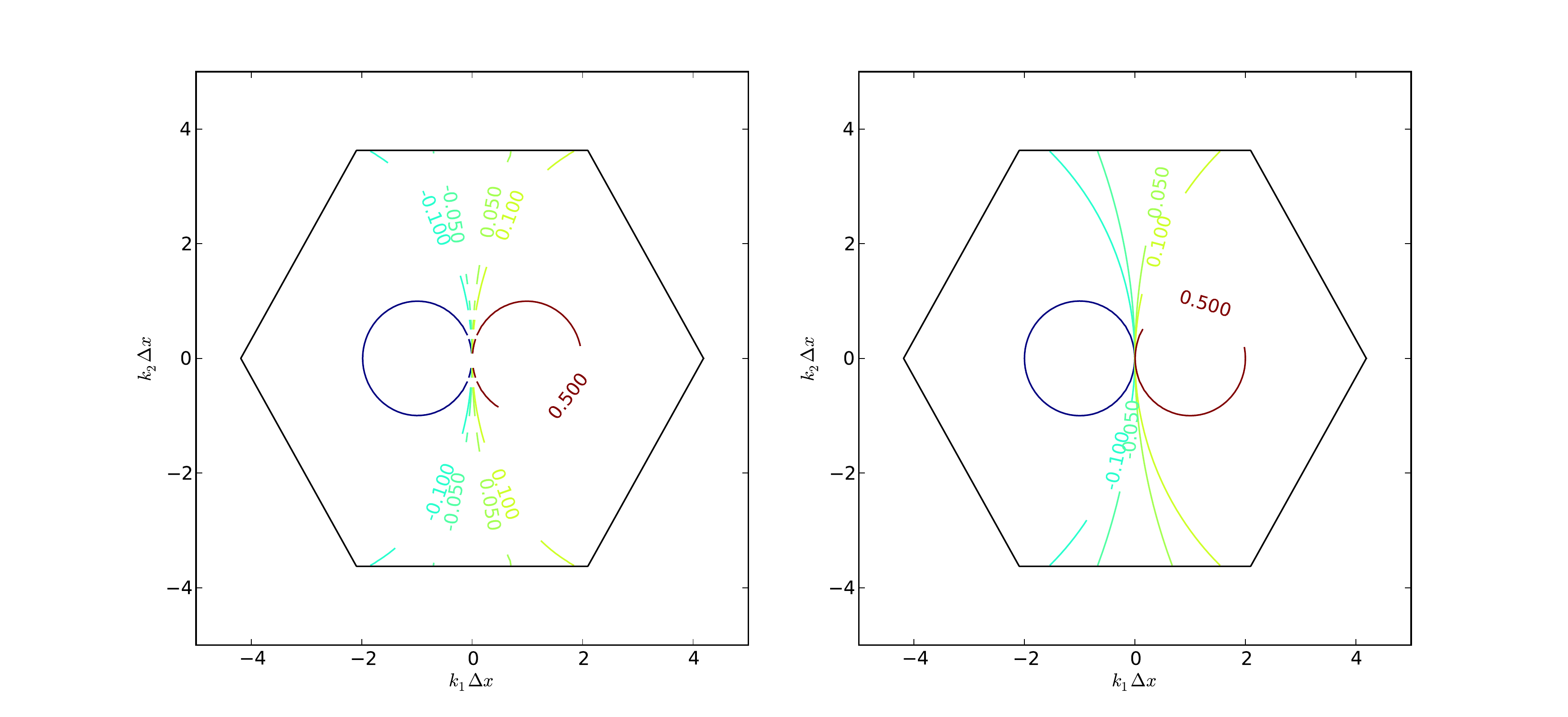}
\end{center}
\caption{\label{rossby_horiz_compare_exact}Comparison between the
  lowest branch of the discrete dispersion relation (left) and the
  exact dispersion relation (right). $f$ increases in the
  $y$-direction relative to the mesh.}
\end{figure}

\begin{figure}
\begin{center}
\includegraphics*[width=10cm]{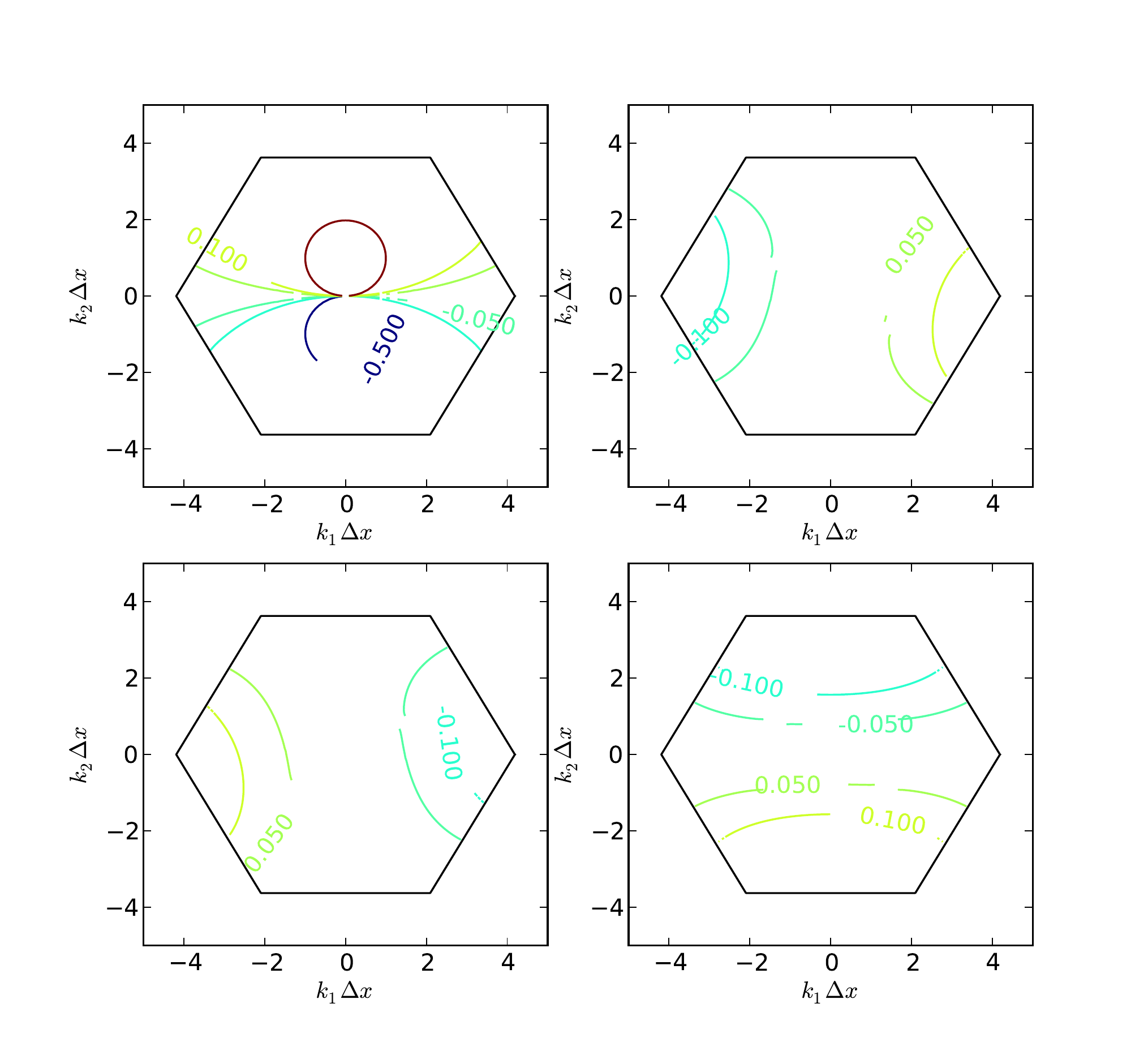}
\end{center}
\caption{\label{rossby_vert_all_modes} Contour plots showing
  $\omega\times 10^6$ obtained from the solutions of equation
  \eqref{rossby dispersion relation}, with parameters $f_0 =
  1.0\times10^{-4}$, $\beta = 1.0\times 10^{-12}$, $\Delta x =
  1.0\times 10^5$ and $c^2=1.0\times 10^5$ (these parameters are the
  same as those used in \citet{Th08}). $f$ increases in the
  $x$-direction relative to the mesh. The lowest branch of the
  dispersion relation is shown top-left. The other branches are
  aliased higher values of $\MM{k}\Delta x$.}
\end{figure}

\begin{figure}
\begin{center}
\includegraphics*[width=10cm]{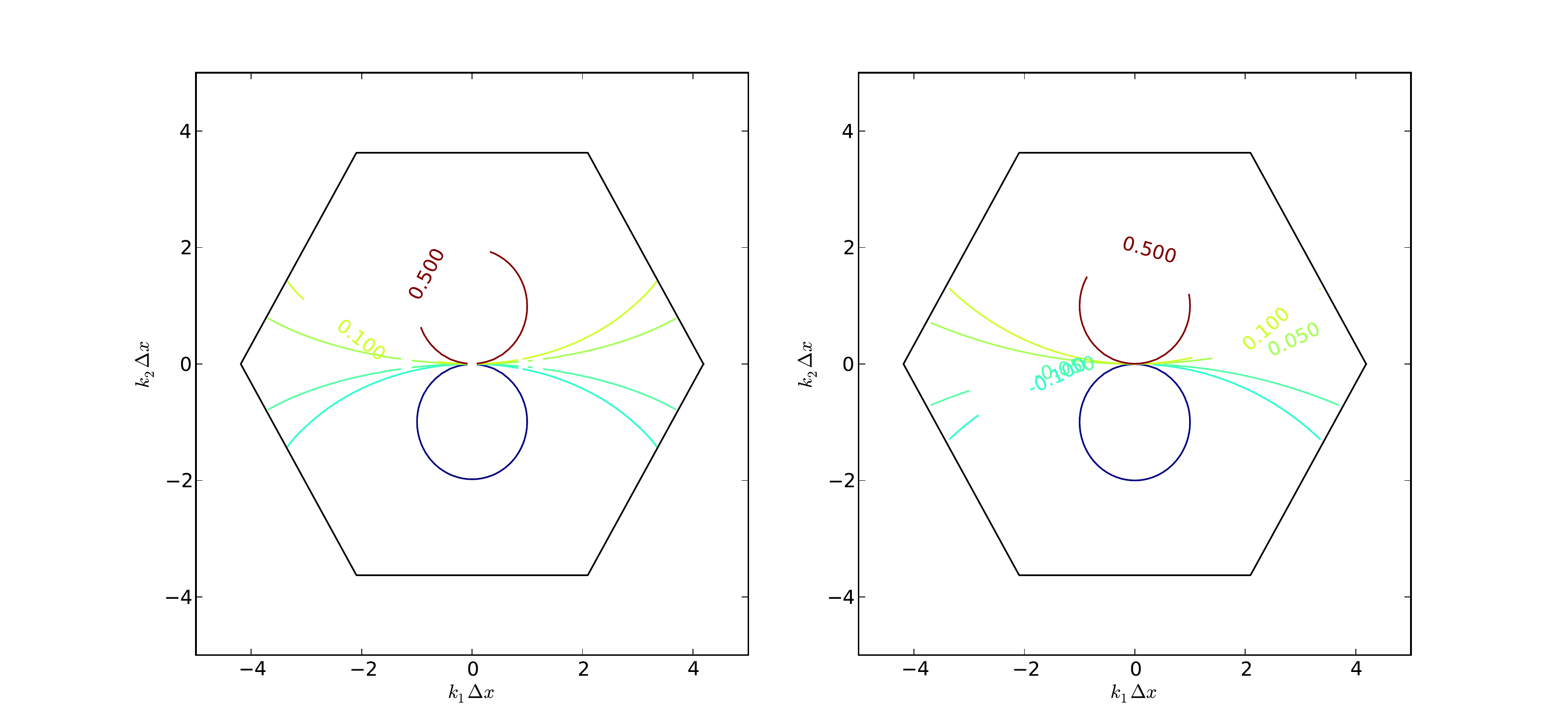}
\end{center}
\caption{\label{rossby_vert_compare_exact}Comparison between the
  lowest branch of the discrete dispersion relation (left) and the
  exact dispersion relation (right). $f$ increases in the
  $y$-direction relative to the mesh.}
\end{figure}

Figure \ref{rossby_horiz_all_modes} shows contour plots of the
frequency $\omega$ for the case $\hat{\MM{f}}=(0,1)$, with parameter
values taken from \citet{Th08}.  Exactly as the $f$-plane case, we
obtain four roots for $\omega$ which correspond to the fundamental
modes (\emph{i.e.} the modes that are possible to represent on a $P1$
mesh) and the higher modes which arise from the extra accuracy on a
$P2$ mesh. All the modes correspond to physical values after correct
interpretation through the Brillouin zones as for the inertia-gravity
wave case. A comparison with the exact dispersion relation for Rossby
waves is given in Figure \ref{rossby_horiz_compare_exact}; a very
close match is observed.  Figure \ref{rossby_vert_all_modes} shows
contour plots for the same parameter values but with
$\hat{\MM{f}}=(1,0)$.  Figure \ref{rossby_vert_compare_exact} shows
the corresponding comparison with the exact dispersion relation; a
close match is again observed. This shows that the \pdgp
discretisation has Rossby waves whose speed is almost independent of
the mesh orientation.

\section{Summary and outlook}
\label{summary}


In this paper we analysed the \pdgp finite element pair applied to the
rotating shallow-water equations, by means of a discrete Helmholtz
decomposition which exists because of the embedding properties of
\pdgp, namely gradients and skew gradients of \ptwo map into \pdg.
The discrete Helmholtz decomposition has some extra components, which
we refer to as spurious velocity components, and which can be
projected out, resulting in a discretisation that we referred to as
\hptwo-\ptwo. This decomposition was then used to show that in the
$f$-plane, all steady states are geostrophically balanced (and
\emph{vice versa}). Furthermore, a discrete inertia-gravity wave
equation can be derived which is the same as the \ptwo continuous
finite element method applied to the inertia-gravity wave equation,
and hence the inertia-gravity wave solutions are third-order
accurate. This should mean that the \pdgp method should give very
stable and accurate solutions of the linear geostrophic adjustment
problem. We also showed that the spurious velocity components are
uncoupled from the geostrophic balance or inertia-gravity waves, and
they just undergo spurious inertial oscillations which do not
propagate. When the \hptwo-\ptwo method is used, we obtain identical
equations but without the spurious inertial oscillations. The
\hptwo-\ptwo method may be thought of as an implementation of the
\ptwo finite element version of the Z-grid, in which vorticity,
streamfunction and layer thickness are all collocated. Hence, the
\pdgp method may be thought of as a way to embed the finite element
Z-grid into a method which avoids the need to solve elliptic problems
for the potential and streamfunction, at the cost of adding spurious
inertial oscillations.

We then followed the methodology of \citet{LePo2008,Th08} to analyse
the Rossby wave equation obtained from the the \pdgp discretisation of
the shallow-water equations on the $\beta$-plane in the
quasi-geostrophic limit. It was shown that the spurious velocity
components do not couple in to the Rossby wave dynamics, in fact the
geostrophic spurious components vanish and the ageostrophic components
are slaved to the geostrophic streamfunction. It was shown that the
quasi-geostrophic limit leads to a discrete Rossby wave equation which
is identical to the continuous \ptwo finite element discretisation
applied to the continuous Rossby wave equation, and hence the \pdgp
Rossby waves are third-order accurate. We expect that this means that
the \pdgp dispersion relation is much more independent of the
direction of mesh alignment than other methods with linear velocity
(such as the lowest-order Ravier-Thomas element which is the finite
element version of the C-grid finite difference method). One seemingly
negative aspect of using continuous finite element methods for
pressure is that the mass matrix is not diagonal, so a linear system
must be solved even when explicit timestepping is used. On the one
hand, solving this linear system iteratively is extremely cheap since
the condition number is independent of resolution and hence the number
of iterations required stays constant under mesh refinement
\citep{GrSa2000}. On the other hand, one can approximate the mass
matrix $M$ by a ``lumped'' diagonal mass matrix $M_L$ with
$(M_L)_{ii}=\sum_{j}M_{ij}$. It was shown in \citep{LeHaRoPo2008} that
lumping the mass has minimal effect on the dispersion relations so we
would expect similar properties. In particular note that mass-lumping
only effects the time-derivative terms so geostrophic states will
remain steady.

It seems almost inevitable (because of the difficulty in balancing the
number of velocity and pressure degrees of freedom) that any numerical
discretisation that is not based on quadrilateral meshes will result
in some form of spurious modes. From the results of this paper it
appears that the \pdgp method puts the spurious modes into the least
harmful place: it has no spurious pressure modes which would quickly
pollute the solution and result in sub-optimal numerical convergence,
it has no spurious Rossby modes which could modify the transfer of
energy from barotropic to baroclinic modes in the presence of
baroclinic instability, but it does have spurious inertial
oscillations which do not propagate, and which can be filtered out
using the \hptwo-\ptwo projection. Whether or not these modes cause
problems depends on how they are coupled to the physical modes through
nonlinear advection, and this needs to be studied in careful
benchmarks before recommending the \pdgp method for use in NWP.  If
the modes are not harmful then the other properties discussed here
(super-accurate wave propagation and representation of geostrophic
balance on arbitrary unstructured meshes) mean that \pdgp should be an
ideal choice for NWP models using adaptive mesh refinement.  Here the
projection filter will prove very useful, since the spurious modes can
easily be extracted and measured, and modified advection schemes can
be proposed which apply the projection before the wave step in
semi-implicit splitting methods.

\paragraph{Acknowledgements} This paper began after interesting
discussions on spurious modes with Andrew Staniforth, John Thuburn and
Nigel Wood.

\bibliographystyle{elsarticle-harv}
\bibliography{wavepropagationP1dgP2}

\begin{thebibliography}{22}
\expandafter\ifx\csname natexlab\endcsname\relax\def\natexlab#1{#1}\fi
\expandafter\ifx\csname url\endcsname\relax
  \def\url#1{\texttt{#1}}\fi
\expandafter\ifx\csname urlprefix\endcsname\relax\def\urlprefix{URL }\fi

\bibitem[{Arakawa and Lamb(1977)}]{ArLa77}
Arakawa, A., Lamb, V., 1977. Computational design of the basic dynamical
  processes of the {UCLA} general circulation model. In: Chang, J. (Ed.),
  Methods in Computational Physics. Vol.~17. Academic Press, pp. 173--265.

\bibitem[{Brenner and Scott(1994)}]{BrSc04}
Brenner, S., Scott, R., 1994. The Mathematical Theory of Finite Element
  Methods. Springer-Verlag.

\bibitem[{Comblen et~al.(2010)Comblen, Lambrechts, Remacle, and
  Legat}]{CoLaReLe2010}
Comblen, R., Lambrechts, J., Remacle, J.-F., Legat, V., 2010. Practical
  evaluation of five partly discontinuous finite element pairs for the
  non-conservative shallow water equations. Int. J. Num. Meth. Fluid. 63~(6),
  701--724.

\bibitem[{Cotter et~al.(2009{\natexlab{a}})Cotter, Ham, and Pain}]{CoHaPa2009}
Cotter, C.~J., Ham, D.~A., Pain, C.~C., 2009{\natexlab{a}}. A mixed
  discontinuous/continuous finite element pair for shallow-water ocean
  modelling. Ocean Modelling 26, 86--90.

\bibitem[{Cotter et~al.(2009{\natexlab{b}})Cotter, Ham, Pain, and
  Reich}]{CoHaPaRe2009}
Cotter, C.~J., Ham, D.~A., Pain, C.~C., Reich, S., 2009{\natexlab{b}}. {LBB}
  stability of a mixed finite element pair for fluid flow simulations. J. Comp.
  Phys. 228~(3), 336--348.

\bibitem[{Fox-Rabinovitz(1996)}]{Fo-Ra96}
Fox-Rabinovitz, M., 1996. Computational dispersion properties of {3D} staggered
  grids for a nonhydrostatic anelastic system. Mon. Weather Rev. 124, 498--510.

\bibitem[{Gresho and Sani(2000)}]{GrSa2000}
Gresho, P.~M., Sani, R.~L., 2000. Incompressible Flow and the Finite Element
  Method, Volume 2, Isothermal Laminar Flow. Wiley.

\bibitem[{Kossevich(2005)}]{Ko2005}
Kossevich, A.~M., 2005. Geometry of Crystal Lattices. Wiley.

\bibitem[{{Le Roux} et~al.(2008){Le Roux}, Hanert, Rostand, and
  Pouliot}]{LeHaRoPo2008}
{Le Roux}, D., Hanert, E., Rostand, V., Pouliot, B., 2008. Impact of mass
  lumping on gravity and rossby waves in 2d finite-element shallow-water
  models. Int. J. Num. Meth. Fluid. 59~(7), 767--790.

\bibitem[{Le~Roux et~al.(1998)Le~Roux, Staniforth, and Lin}]{Ro_etal1998}
Le~Roux, D., Staniforth, A., Lin, C.~A., 1998. Finite elements for
  shallow-water equation ocean models. Monthly Weather Review 126~(7),
  1931--1951.

\bibitem[{Majewski et~al.(2002)Majewski, Liermann, Prohl, Ritter, Buchhold,
  Hanisch, Paul, Wergen, and Baumgardner}]{Ma+2002}
Majewski, D., Liermann, D., Prohl, P., Ritter, B., Buchhold, M., Hanisch, T.,
  Paul, G., Wergen, W., Baumgardner, J., 2002. The operational global
  icosahedral-hexagonal gridpoint model {GME}: Description and high-resolution
  tests. Mon. Wea. Rev. 130, 319--338.

\bibitem[{Randall(1994)}]{Ra94}
Randall, D., 1994. Geostrophic adjustment and the finite-difference
  shallow-water equations. Mon. Weather Rev. 122, 1371--1377.

\bibitem[{Raviart and Thomas(1977)}]{RaTh1977}
Raviart, Thomas, 1977. A mixed finite element method for 2nd order elliptic
  problems. In: Mathematical Aspects of the Finite Element Method. Lecture
  Notes in Mathematics. Springer, Berlin, pp. 292--315.

\bibitem[{Ringler et~al.(2000)Ringler, Heikes, and Randall}]{Ri+2000}
Ringler, T.~D., Heikes, R., Randall, D., 2000. Modeling the atmospheric general
  circulation using a spherical geodesic grid: A new class of dynamical cores.
  Mon. Wea. Rev. 128, 2471--2490.

\bibitem[{Roux and Pouliot(2008)}]{LePo2008}
Roux, D. Y.~L., Pouliot, B., 2008. Analysis of numerically induced oscillations
  in two-dimensional finite-element shallow-water models part ii: Free
  planetary waves. SIAM Journal on Scientific Computing 30~(4), 1971--1991.

\bibitem[{Roux et~al.(2007)Roux, Rostand, and Pouliot}]{LeRoPo2007}
Roux, D. Y.~L., Rostand, V., Pouliot, B., 2007. Analysis of numerically induced
  oscillations in 2d finite-element shallow-water models part i:
  Inertia-gravity waves. SIAM Journal on Scientific Computing 29~(1), 331--360.

\bibitem[{Satoh et~al.(2008)Satoh, Matsuno, Tomita, Miura, Nasuno, and
  Iga}]{Sa+2008}
Satoh, M., Matsuno, T., Tomita, H., Miura, H., Nasuno, T., Iga, S., 2008.
  Nonhydrostatic icosahedral atmospheric model {(NICAM)} for global cloud
  resolving simulations. J. Comp. Phys. 227~(7), 3486--3514.

\bibitem[{{SymPy Development Team}(2009)}]{SymPy}
{SymPy Development Team}, 2009. SymPy: Python library for symbolic mathematics.
\newline\urlprefix\url{http://www.sympy.org}

\bibitem[{Thuburn(2008)}]{Th08}
Thuburn, J., 2008. Numerical wave propagation on the hexagonal {C}-grid. J.
  Comp. Phys. 227~(11), 5836--5858.

\bibitem[{Thuburn et~al.(2009)Thuburn, Ringler, Skamarock, and
  Klemp}]{ThRiSkLe2009}
Thuburn, J., Ringler, T.~D., Skamarock, W.~C., Klemp, J.~B., 2009. Numerical
  representation of geostrophic modes on arbitrarily structured {C}-grids. J.
  Comput. Phys. 228, 8321--8335.

\bibitem[{Umgiesser et~al.(2004)Umgiesser, Canu, Cucco, and Solidoro}]{Um+2004}
Umgiesser, G., Canu, D.~M., Cucco, A., Solidoro, C., 2004. A finite element
  model for the {V}enice {L}agoon. {D}evelopment, set up, calibration and
  validation. Journal of Marine Systems 51~(1-4), 123--145.

\bibitem[{Walters and Casulli(1998)}]{WaCa1998}
Walters, R., Casulli, V., 1998. A robust, finite element model for hydrostatic
  surface water flows. Communications in Numerical Methods in Engineering 14,
  931--940.

\end{thebibliography}

\end{document}